\documentclass[amssymb,amsfonts,12pt]{amsart}
\newtheorem{thm}{Theorem}[section]
\newtheorem{lemma}[thm]{Lemma}

\newtheorem{prop}[thm]{Proposition}
\newtheorem{definition}[thm]{Definition}

\theoremstyle{remark}
\newtheorem{remark}[thm]{Remark}

\usepackage{amsmath}
\usepackage{amsfonts}
\usepackage{mathrsfs}
\usepackage{amssymb}
\usepackage{amsthm}
\usepackage{graphicx}
\usepackage{epic}
\usepackage{eepic}

\usepackage{color}
\usepackage{verbatim}
\usepackage{fixmath}
\usepackage{pstricks}
\usepackage{multido}

\def\QSet{\mbox{\rm\kern.24em
\vrule width.03em height1.48ex depth-.051ex \kern-.26em Q}}

\def\R{{\mathbb R}}

\def\N{{\mathbb N}}
\def\C{{\mathbb C}}

\def\Z{{\mathbb Z}}
\def\l{{\mathcal L}}

\def\be#1{\begin{equation}\label{#1}}

\def\bas{\begin{align*}}
\def\eas{\end{align*}}
\def\bi{\begin{itemize}}
\def\ei{\end{itemize}}

\def\emph#1{{\it #1}}

\newcommand{\SSS}{\mathbb{S}}

\newcommand{\lb}{\langle}
\newcommand{\rb}{\rangle}

\title[Vector valued inequalities]{Vector valued inequalities for families of bilinear Hilbert transforms and applications to bi-parameter problems}
\author{Prabath Silva}\address{Department of Mathematics, Indiana University, 831 East 3rd St., Bloomington IN 47405}
\email{pssilva@indiana.edu}
\begin{document}
\begin{abstract}
Muscalu, Pipher, Tao and Thiele \cite{MPTT} showed that the tensor product between two one dimensional paraproducts (also known as bi-parameter paraproduct) satisfies all the expected $L^p$ bounds. In the same paper they showed that the tensor product between two bilinear Hilbert transforms is unbounded in any range. They also raised the question about $L^p$ boundedness of the bilinear Hilbert transform tensor product with a paraproduct. We answer their question by obtaining a wide range of estimates for this hybrid bilinear operator. Our method relies on new vector valued estimates for a family of bilinear Hilbert transforms.
\end{abstract}

\maketitle

\section{Introduction}

Let $\Gamma_1$ and $  \Gamma_2$ be linear subspaces of $\R^2$, and  consider the smooth symbol  $m:\R^4\to\C$ which satisfies  the  condition
\begin{equation}\label{HMcondition}
|\partial^{\alpha}_{\xi} \partial^{\beta}_{\eta} m (\xi, \eta)| \lesssim \frac{1}{\text{dist}(\xi, \Gamma_1)^{\alpha}}  \frac{1}{\text{dist}(\eta, \Gamma_2)^{\beta}}
\end{equation}
for $ \xi \notin \Gamma_1, \eta \not \in \Gamma_2$.
Consider the  bilinear operator   $T_m$  associated to $m$ defined by
\begin{equation}\label{defTm}
T_m(f,g)(x,y)= \int m(\xi, \eta)  \hat{f}(\xi_1,\eta_1)\hat{g}(\xi_2,\eta_2) e^{2\pi i (x,y) \cdot ((\xi_1,\eta_1) + (\xi_2,\eta_2))} d\xi d\eta,
\end{equation}
where $\hat{f}$ denotes the Fourier transform of $f$,
$$ \hat{f}(\xi,\eta)= \int_{\R^2} f(x,y) e^{-i(x,y)\cdot (\xi,\eta)} dxdy.$$

When  $\dim(\Gamma_1)= \dim(\Gamma_2)=0$ the operator $T_m$ is called a bi-parameter paraproduct and the following theorem is proved in \cite{MPTT}:

\begin{thm}[Bi-parameter Paraproduct]
\label{vjkhrtug5igujkuy}
Consider $m$ and $T_m$ defined as above with $\dim(\Gamma_1)= \dim(\Gamma_2)=0$. Then $T_m$ maps $L^p \times L^q$ into $L^r$ when $\frac{1}{p}+\frac{1}{q}=\frac{1}{r}$ and $p,q> 1.$
\end{thm}

The double bilinear Hilbert transform corresponds to the case when $\dim(\Gamma_1)= \dim(\Gamma_2)=1$. A prototypical example is when $m(\xi,\eta)=sgn(\xi_1-\xi_2)sgn(\eta_1-\eta_2)$, in which case there is a nice kernel representation
$$T_m(f,g)(x,y)=\int\int_{\R^2}f(x+t,y+s)g(x-t,y-s)\frac{dt}{t}\frac{ds}{s}.$$

It was proved in \cite{MPTT} that this operator does not satisfy any $L^p$ estimates. This is a bit surprising, since this operator can be thought of as being a tensor product -in a bilinear setting- of two  bilinear Hilbert transforms
$$B(f,g)(z)=\int_{\R}f(z+t)g(z-t)\frac{dt}{t}.$$
We recall the celebrated theorem  proved in \cite{LT1}.
\begin{thm}\label{BHTthm}
The bilinear Hilbert transform maps $L^p \times L^q$ into $L^r$ when $\frac{1}{p}+\frac{1}{q}=\frac{1}{r}$ with $p,q> 1,  r>\frac{2}{3}.$
\end{thm}

In this paper we consider the case  $\dim(\Gamma_1)=1, \dim(\Gamma_2)=0$. We call $\Gamma_1$  non degenerate when $\Gamma_1 \not = \R \times \{0 \}, \{0\} \times \R$; the analysis in the degenerate cases is easier and we omit it. In the non degenerate case  we can view $T_m$  as a tensor product between the bilinear Hilbert transform and a  paraproduct. We will denote this operator by $BP$ for the rest of the paper.

We prove the following $L^p$ estimates for $BP$. This answers Question 8.2 in \cite{MPTT}.

\begin{thm}[Main Theorem]\label{Main-Theorem}
The operator BP maps $L^p \times L^q$ into $L^r$ with $\frac{1}{p}+\frac{1}{q}=\frac{1}{r}$ and $p,q> 1, \frac{1}{p} +\frac{2}{q} <2, \frac{1}{q} +\frac{2}{p} <2.$
\end{thm}

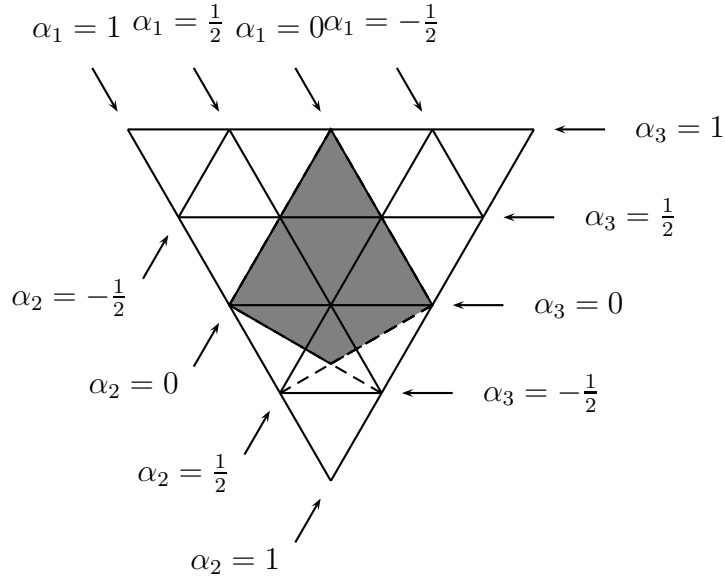
\begin{figure}[h]
\psset{unit=1.35cm}
\begin{pspicture}(-2,1.7321)(6,-5.1961)
\pspolygon[linestyle=dashed,fillcolor=gray,fillstyle=solid](1,-1.732)(2,0)(3,-1.732)(2, -2.3093)
\multido{\n=0+.5,\nn=4+-.5,\nnn=0+-.8660}{4}{\psline(\n,\nnn)(\nn,\nnn)}
\multido{\n=0+1,\nn=2+.5,\nnn=-3.4641+.8660}{4}{\psline(\n,0)(\nn,\nnn)}
\multido{\n=4+-1,\nn=2+-.5,\nnn=-3.4641+.8660}{4}{\psline(\n,0)(\nn,\nnn)}
\multido{\n=4.2+-.5,\nn=4.7+-.5,\nnn=0+-.8660}{4}{\psline{<-}(\n,\nnn)(\nn,\nnn)}
\pspolygon[linestyle=dashed](1,-1.732)(2,0)(3,-1.732)(2, -2.3093)
\psline[linestyle=dashed](1,-1.732)(2.5,-2.6)
\psline[linestyle=dashed](3,-1.732)(1.5,-2.6)

\multido{\n=-.35+1.00,\nn=-.1+1.0}{4}{\psline{->}(\n,.6062)(\nn,.1732)}
\multido{\n=.15+.50,\nn=.4+.5,\nnn=-1.4722+-.8660,\nnnn=-1.0392+-.8660}{4}{\psline{->}(\n,\nnn)(\nn,\nnnn)}
\rput[l](5,0){$\alpha_3=1$}
\rput[l](4.5,-.8660){$\alpha_3=\frac{1}{2}$}
\rput[l](4,-1.732){$\alpha_3=0$}
\rput[l](3.5,-2.5981){$\alpha_3=-\frac{1}{2}$}
\rput[Br](0,-1.732){$\alpha_2=-\frac{1}{2}$}
\rput[Br](.5,-2.5981){$\alpha_2=0$}
\rput[Br](1,-3.4641){$\alpha_2=\frac{1}{2}$}
\rput[Br](1.5,-4.3301){$\alpha_2=1$}
\rput[b](-.5,.8660){$\alpha_1=1$}
\rput[b](.5,.8660){$\alpha_1=\frac{1}{2}$}
\rput[b](1.5,.8660){$\alpha_1=0$}
\rput[b](2.5,.8660){$\alpha_1=-\frac{1}{2}$}
\end{pspicture}\caption{Range of exponents for T.}\label{exponentsT}
\end{figure}

\begin{remark}
Figure \ref{exponentsT} shows the full range of exponents in the plane $\alpha_1+\alpha_2+\alpha_3=1$; there $\alpha_1 =\frac{1}{p}, \alpha_2 =\frac{1}{q}$, $ \alpha_3 = \frac{1}{r'}$ while $r'$ denotes the dual exponent of $r$. Note that this range of exponents includes the Banach triangle, that is the range $p,q,r>1$.
\end{remark}

We prove Theorem \ref{Main-Theorem} in Section \ref{proofofmainthm} . To illustrate the main circle of ideas behind the proof of Theorem \ref{Main-Theorem}, we first prove the less technically involved local $L^2$ case  in sections \ref{l2} and \ref{VVBHTPsec}.

We would like to point out the fact that Theorem \ref{Main-Theorem} has resisted many attempts by various mathematicians to employ time-frequency techniques with product theory flavor. These techniques include product trees, product BMO, product John-Nirenberg and Journe's lemma. While the approach in \cite{MPTT} relies, at the very minimum,  on the boundedness of the square function on the bi-disk, our whole approach in this paper only assumes the boundedness of the strong maximal function. Thus our methods are very one dimensional in nature.

The proof of Theorem \ref{Main-Theorem}  uses  vector valued estimates for the  single scale operators $BP_j$ defined below.

We start by recalling the Littlewood-Paley decomposition. Let $\phi$ be a smooth bump function with $ 1_{(-1, 1)} \le \phi \le 1_{(-\frac{21}{20}, \frac{21}{20})}$ and let $ \psi(x) = \phi(x/2) -\phi (x)$. Now we have
$$1_{\R} =\phi + \sum_{j \ge 1}\psi_j$$
$$1_{\R \setminus \{0 \} } = \sum_{j \in \Z}\psi_j $$
where $\phi_j$ and $\psi_j$ are $L^{\infty}$ dilations defined by $\phi_j(x) =\phi(\frac{x}{2^j})$ and $\psi_j(x) =\psi(\frac{x}{2^j})$. Note also that $\text{supp } \psi_j \subset [-2^{j+2},- 2^j]\cup[2^j, 2^{j+2}]$.

\begin{definition}\label{LPproj}
We define the Littlewood-Paley frequency projections at scale $2^j$ in the second coordinate by
$$\Pi^l_j f = \left( \hat{f} (\xi, \eta ) \psi_j(\eta)\right)^\vee$$
$$\Pi^o_j f = \left( \hat{f} (\xi, \eta) \phi_j(\eta)\right)^\vee.$$

Define also the single scale operators
\begin{equation}
\label{defTj}
BP_j(f,g) =BP_j(\Pi^o_j f, \Pi^l_j g).
\end{equation}
\end{definition}

The following vector valued estimates proved in Section \ref{VVBHTPsec} are the key ingredient in dealing with the local $L^2$ case of Theorem \ref{Main-Theorem}.

\begin{thm}\label{VVBHTPl2} Let $f= \{f_j \}$, $g =\{g_j\}$ be sequences of complex valued functions on $\R^2$. Then we have the following estimate for $2<p,q,r'< \infty$,
\begin{equation}
\|\left\{ BP_j(f_j,g_j) \right\}_j\|_{L^r(l^2)} \lesssim  \left\| f  \right\|_{L^p(l^{\infty})}\left\|  g \right\|_{L^q(l^2)}
\end{equation}
\end{thm}

When using $BP_j$ to get estimates for $BP$ it is important to note that $BP_j$ is not symmetric with respect to the entries $f$ and $g$, in fact this is one of the main issues that we have to deal with when we prove the main Theorem \ref{Main-Theorem}. Also, when we apply Theorem \ref{VVBHTPl2}  to get estimates for $BP$ we take $f_j, g_j$ to be Littlewood-Paley frequency projections, and so it is important to work with $l^2$ and $l^\infty$. These are in fact the only two $l^p$ spaces that will have a relevance for proving the main Theorem \ref{Main-Theorem}.

However, our methods can prove estimates in a larger range as we see below, see also Figures \ref{ARR=2}, \ref{ARR<2} and \ref{ARR>2}.  First for $\frac{4}{3} < R < 4$ define,
\begin{equation}
A^R =\left\{\alpha \in \pi_0: 0< \alpha_1<1-4\left| \frac{1}{2}- \frac{1}{R}\right| , |\alpha_2-\alpha_3|<1-2\left| \frac{1}{2}- \frac{1}{R}\right|  \right\}.
\end{equation}
We refer the reader to Definition \ref{exp} for $\pi_0$.
\begin{thm}\label{VVBHTPmess}
We have
\begin{equation}\label{VVBHTPmessest}
\|\left\{ BP_j(f_j,g_j)_j \right\}\|_{L^r(l^R)} \lesssim  \left\| f  \right\|_{L^p(l^{\infty})}\left\|  g \right\|_{L^q(l^R)}
\end{equation}
when $\frac{1}{p}+\frac{1}{q}=\frac{1}{r}$, $p,q >1$, $\frac{4}{3} < R < 4$, and $(\frac{1}{p}, \frac{1}{q}, \frac{1}{r'}) \in A^R. $
 \end{thm}

We prove an interpolation result in Section \ref{grtLPsec}, which allow us to interpolate inner spaces $l^P$ using complex interpolation and outer spaces $L^p$ using generalized restricted type interpolation. We use this interpolation result to get estimates in some even larger range.

To see this extension, let us introduce some notation. For $\alpha=(\alpha_1, \alpha_2, \alpha_3)\in \pi_0$, we define $\bar{\alpha} =(\alpha_2, \alpha_1, \alpha_3)$. Next for given $\frac43<R<4$ define $\mathbf{A}^R \subset \pi_0 \times \pi_0$ by,
$$\mathbf{A}^R = \{ (\boldsymbol{\alpha}, \alpha ): \boldsymbol{\alpha}=(0, \frac{1}{R}, \frac{1}{R'}), \alpha \in A_R \}.$$
We define $\mathcal{A}^R \subset \pi_0 \times \pi_0$ by
\begin{equation}\label{range}
\mathcal{A}^R = \{ (\boldsymbol{\alpha}, \alpha ): (\boldsymbol{\alpha}, \alpha) \in \mathbf{A}^R \text{ or }   (\bar{\boldsymbol{\alpha}}, \bar{\alpha}) \in \mathbf{A}_R \}
\end{equation}
and let $Converx(B)$ denote the convex hull of a set $B$.
\begin{thm}\label{VVBHTPmess1}
We have
\begin{equation}\label{VVBHTPest}
\left\| \left( \sum |BP_j(f_j,g_j)|^R \right)^{\frac{1}{R}} \right\|_r \lesssim \left\| \left( \sum |f_j|^P \right)^{\frac{1}{P}} \right\|_p\left\| \left( \sum |g_j|^Q \right)^{\frac{1}{Q}} \right\|_q
\end{equation}
for $\frac{1}{p}+\frac{1}{q}=\frac{1}{r}$ with $p,q> 1$ and $\frac{1}{P}+\frac{1}{Q} = \frac{1}{R}$, $P,Q>1$, $\frac{4}{3} < R < 4$ and $(\frac{1}{P}, \frac{1}{Q}, \frac{1}{R'},   \frac{1}{p}, \frac{1}{q}, \frac{1}{r'}) = (\boldsymbol{\alpha}, \alpha ) \in Convex(\mathcal{A}^R)$.
 \end{thm}

Our methods reduce the vector valued inequalities to estimates for single operators $BP_j$. We use the fact that each operator $BP_j$ satisfies certain time frequency estimates outside appropriate exceptional sets. In these regards, our approach bears some resemblance to the argument in \cite{BT1}, though the bilinear nature of our estimates creates new difficulties. One of these difficulties is encountered in the proof of estimates outside the Banach triangle, in Section \ref{proofofmainthm}. The fact that one can have estimates below $L^1$ is a feature unique to the multilinear setting. The unboundedness of the Littlewood-Paley square function below $L^1$ calls for extra care in the argument.

We close this preliminary discussion by putting our results in a more general framework. Let us consider a family of kernels $K_j$ which satisfy the decay estimates,
\begin{equation}\label{Kj}
|\partial^{\alpha}\hat{K_j}(\xi)| \le \frac{C_{\alpha}}{|\xi|^{\alpha}}
\end{equation}
for $\alpha \in \N,$ where $C_{\alpha}$ is independent of $j$. Define the bilinear operator  $B_j$ associated to the kernel  $K_j$ by
\begin{equation}\label{BHTKj}
B_j(f,g)(x) = \int_{\R} f(x-t)g(x+t)K_j(t)dt.
\end{equation}
A particular example is $K(t)=1/t$. And indeed, all these operators  fall under the scope of Theorem \ref{BHTthm}, and satisfy the same estimates, uniformly in $j$.

We have the following vector valued version of Theorem \ref{BHTthm} for the family $B_j$ of bilinear Hilbert transforms.
\begin{thm}\label{VVBHT}
We have
\begin{equation}\label{VVBHTest}
\left\| \left( \sum |B_j(f_j,g_j)|^R \right)^{\frac{1}{R}} \right\|_r \lesssim \left\| \left( \sum |f_j|^P \right)^{\frac{1}{P}} \right\|_p\left\| \left( \sum |g_j|^Q \right)^{\frac{1}{Q}} \right\|_q
\end{equation}
for $\frac{1}{p}+\frac{1}{q}=\frac{1}{r}$ with $p,q> 1$ and $\frac{1}{P}+\frac{1}{Q} = \frac{1}{R}$, $P,Q>1$, $\frac{4}{3} < R < 4$ and $(\frac{1}{P}, \frac{1}{Q}, \frac{1}{R'},   \frac{1}{p}, \frac{1}{q}, \frac{1}{r'}) = (\boldsymbol{\alpha}, \alpha ) \in \text{Convex} (\mathcal{A}^R)$.
 \end{thm}

Let us mention that in the one kernel case $K_j=K$, the above estimate  when $P=Q=R=2$, $p,q>1$ and  $r>\frac23$, follows immediately from Theorem \ref{BHTthm}, using the standard randomization argument. This has been observed in \cite{GL}. We point out again that in order to derive our main application, Theorem \ref{Main-Theorem}, we need to work with families of kernels $K_j$, rather than just one kernel. This necessity comes in part from the presence of the coefficients $\epsilon_s$ in \eqref{eq: modelsum}. Equally importantly, due to the presence of one overlapping index, working with just $l^2$ is not enough for our purposes. We need to bring the $l^\infty$ space into the picture.

By using the methods in this paper, we can recover part of the range from \cite{MPTT}, for the bi-parameter paraproduct. More precisely, we can show that Theorem \ref{vjkhrtug5igujkuy} holds true under the additional restriction $r > \frac32$.

An interesting question is whether the range in Theorem \ref{Main-Theorem} can be pushed to match that from Theorem \ref{BHTthm}. While this is somewhat natural to expect, our methods do not seem to shed light on this issue.

\subsection{Acknowledgements:} The author would like to thank his thesis advisor, Ciprian Demeter, for his guidance and many helpful conversations about this problem.

\section{Notation}

In this section we introduce some notation that we use throughout the paper. We use $A \lesssim B$ to denote $A\le C B$, where $C$ is a large absolute constant.

Given a cube $Q \subset \R^n$ (or a rectangle $R \subset \R^2$) we denote the measure of $Q$ (or $R$) by  $|Q|$ (or $|R|$). The side length of $Q$ is denoted by $l(Q)$,  and we use $cQ$ to denote the cube with the same center as $Q$ and $c$ times the side length.

Given $I \subset \R,$ we define the cutoff function $\tilde{\chi_I}$ by
\begin{equation}\label{cutoff}
 \tilde{\chi_I} (x) = (1+ (\frac{|x-x_I|}{|I|})^2)^{-1/2} ,
\end{equation}
where $x_I$ is the center of $I$.

We use $M^1f$ to denote the Hardy-Littlewood maximal function of $f$ on $\R$, and use $Mf$ to denote the strong maximal function on $\R^2$. We use $M_r$ for $r\ge1$ to denote the maximal function defined by
$$M_rf = \left(M|f|^r\right)^{1/r}.$$Let $c(p, q, r)$ be a small but positive absolute constant depending only on $(p, q, r)$.
We use $\epsilon$ to denote a small number dependent  on $p,q,r$ such that $\epsilon > c(p, q, r)$. We caution the reader that our $\epsilon$ can change from line to line in our argument.

\section{Interpolation}\label{BTintersec}

In this section we go over the real interpolation theorems which are used in the paper. From now, $x$ may refer to an element of either $\R$ or $\R^2$. First, for a bilinear operator $T$ define  {\it {the associated trilinear form $ \Lambda =\Lambda_T$ of $T$}}  by
$$ \Lambda(f, g, h) =\int T(f, g)(x) h(x) dx.$$
We recall the generalized restricted type interpolation results for bilinear operators  \cite{MTT1}, (see also \cite{Th1}). These reduce the proof of $L^p$ estimates for bilinear operators to proving the so called {\em generalized restricted  type} estimates for the associated trilinear forms. Then we prove a bilinear variant of the interpolation used in \cite{BT1}, which reduces vector valued estimates for a family of bilinear operators to a family of estimates for the operators in the family.

\begin{definition}\label{exp} A triple $\alpha = (\alpha_1, \alpha_2, \alpha_3)$ is called admissible, if
$$-\infty < \alpha_i <1$$
for all $1\le i\le3$,
$$ \alpha_1+\alpha_2+\alpha_3=1$$
and there is at most one index $j$ such that $\alpha_j<0$. We call an index good if $\alpha_i >0$, and call it bad if $\alpha_i \le 0$. We denote the set of admissible triples by $\pi_0$.

We denote the set of admissible triples without a bad index by $\pi_1$, and we also call such triples  local $L^1$ triples.

The  set of admissible triples with $$ 0< \alpha_i< 1/2$$ for all $1\le i\le3$ is denoted by $\pi_2$, and we also call such triples  local $L^2$ triples.

\end{definition}

\begin{definition} If $E$ is a subset of finite measure of either $\R$ or $\R^2$, we denote by $X(E)$ the space of all measurable complex valued functions $f$ on $\R^2$ (or $\R$) such that $|f| \le 1_E$.
\end{definition}

\begin{definition}
Let $E, E'$ be sets with finite measure. We say that $E'$ is a major subset of $E$ if $E' \subset E$ and $|E' |\ge \frac{1}{2} |E|$. Given a triple $E=(E_1, E_2, E_3)$ and $\alpha \in \pi_0$, we say { \it{ $E'=(E'_1, E'_2, E'_3) $ is a major triple of $(E, \alpha) $}} if there is an index $j \in \{1,2,3 \}$ such that for $i\not =j$ we have $ E'_i =E_i$, $E'_j$ is a major subset of $E_j$, and if $\alpha \in \pi_0 \setminus \pi_1$, then $j$ is the bad index of $\alpha$. Also we refer to $E_j \setminus E_j'$ as an exceptional set of $E_j$.

\end{definition}

\begin{definition}
If $\alpha \in \pi_0$ we say that a trilinear form $\Lambda$ is of  generalized restricted type $\alpha$ with constant $K$ if for triple $E= (E_1,E_2,E_3)$ there exists a major triple $E'$ of $(E, \alpha)$ such that
\begin{equation}\label{grt}
|\Lambda(f_1, f_2, f_3)|  \le K |E|^{\alpha}
\end{equation}
for all functions $ f_i \in X(E'_i), i = 1,2,3.$ Here and in the rest of the paper $|E|^{\alpha} := |E_1|^{\alpha_1}|E_2|^{\alpha_2}|E_3|^{\alpha_3}$.

If  $\Lambda$ is of  generalized restricted type $\alpha$ for  $\alpha \in \pi_1$, we also say $\Lambda$ is of restricted type $\alpha$.
\end{definition}

We record the following lemma (\cite{Th1}, Lemma $3.5$), whose proof follows immediately by simply taking a geometric mean.
\begin{lemma}[Convexity]\label{convexity}
Let $A \subset \pi_0$.
If $\Lambda$ is generalized restricted type $\alpha$ for $\alpha \in A $ with constant $K$,  and the choice of index $j$ and the major subset $E_j'$ is independent of $\alpha$, then $\Lambda$ is generalized restricted type $\beta$ for $\beta \in Convex(A)$ with  the same constant $K$.
\end{lemma}
\begin{remark}\label{noexcinrt}
If $\Lambda$ is  restricted type $\alpha \in \pi_1$ with constant $K$, then by using induction we can get \eqref{grt} for  $E'=E$ with constant $C_{\alpha} K$. We will typically use this fact together with Lemma \ref{convexity},  when $A= \{ \alpha , \alpha_0 \}$ and $\alpha_0 \in \pi_1, \alpha \in \pi_0 \setminus \pi_1$.
\end{remark}

We first state the restricted type interpolation theorem. In our proofs interpolation constants play an important role so we state interpolation theorems with the constants.

\begin{thm}\label{rtthm} Let $\alpha \in \pi_1$ and let $T$ be a bilinear operator. If the associated trilinear form $\Lambda_T$ is  restricted type $\beta \in \pi_1$ in  a neighborhood of $\alpha$ with constant $K,$ then we have
$$\|T(f,g)\|_{r} \lesssim K \|f\|_{p} \|g\|_{q};$$
here $(\frac{1}{p}, \frac{1}{q}, \frac{1}{r'} )= \alpha$ and $\frac{1}{r'}+ \frac{1}{r}=1$.
\end{thm}

Next we state the generalized restricted type interpolation theorem from \cite{MTT1}, see also Theorem 3.8 in \cite{Th1}.

\begin{thm}\label{grtthm} Let $\alpha \in \pi_0\setminus\pi_1$ with first two indices good and let $T$ be a bilinear operator  such that the associated trilinear form $\Lambda_T$ is generalized  restricted type $\beta \in \pi_0$ in  a neighborhood of $\alpha$ with constant $K$.  Then we have
$$\|T(f,g)\|_{r} \lesssim K \|f\|_{p} \|g\|_{q};$$
here $(\frac{1}{p}, \frac{1}{q}, \frac{1}{r'} )= \alpha$ and $\frac{1}{r'}+ \frac{1}{r}=1$.
\end{thm}

Note that according to Lemma \ref{convexity}, to get (generalized) restricted type in a neighborhood of a point $\alpha$, it suffices to check this property for only three points whose open convex hull contains $\alpha$. We will use this observation a few times in our arguments.

Theorems \ref{grtthm} and \ref{rtthm} are about complex valued functions. It turns out that the same proofs will apply, with no essential modifications, to the case when $\C$ is replaced by general Banach spaces. We rewrite  Theorem \ref{grtthm} and Theorem \ref{rtthm} in the case these Banach spaces are $l^P$ spaces, which is the case of interest for our applications. See Lemma \ref{vvrt} and \ref{vvgrt} below.
\begin{definition}\label{XP} For $0 < P \le \infty$ and for a subset  $E$  of $\R^n$ with  finite measure, we denote by $X_{P}(E)$ the space of all sequences of complex valued functions $f = \{ f_j \}$ on $\R^n$ such that $ \|f(x)\|_{l^P} = \left( \sum |f_j(x)|^P \right)^{\frac{1}{P}} \in X(E)$.
\end{definition}

\begin{lemma}\label{vvrt}
Let $\alpha \in \pi_1$ and let $T= \{T_j\}$ be a sequence of bilinear operators. If the following  restricted type condition holds for $\beta \in \pi_1$ in a neighborhood of $\alpha$ with constant $K$ independent of $\beta$---
\begin{quotation}
 For every triple $E=(E_1,E_2,E_3)$ there exists a major triple $E'=(E_1', E_2', E_3')$ of $(E, \beta)$ such that
\begin{equation}\label{VVrtest}
|\int  \left(\sum |T_j(f_j,g_j)|^2 \right)^{\frac{1}{2}} (x)h(x) dx | \le K |E|^{\beta}
\end{equation}
for all sequences of functions  $ f = \{ f_{j}\} \in X_{\infty}(E'_1), g = \{ g_{j}\} \in X_{2}(E'_2), $ and $h \in X(E'_3)$---
\end{quotation}
then we have
$$\| T(f,g) \|_{ L^r(l^2) }:=\left\| \|\left\{ T_j(f_j,g_j) \right\}\|_{l^{2}} \right\|_r \lesssim K \left\| f  \right\|_{L^p(l^\infty)}\left\|  g \right\|_{L^q(l^2)},$$
where $(\frac{1}{p}, \frac{1}{q}, \frac{1}{r'} )= \alpha$ and $\frac{1}{r'}+ \frac{1}{r}=1.$
\end{lemma}

Note that Remark \ref{noexcinrt} applies here too.

\begin{lemma}\label{vvgrt}
Let $\alpha \in \pi_0\setminus\pi_1$ with first two  indices good, and let $T= \{T_j\}$ be a sequence of bilinear operators, which is restricted type $\alpha_0$ for some $\alpha_0 \in \pi_1$ as in Lemma \ref{vvrt}. If the following generalized restricted type condition holds for $\beta \in \pi_0$ in a neighborhood of $\alpha$ with constant $K$ independent of $\beta$---
\begin{quotation}
 For every triple $E=(E_1,E_2,E_3)$ there exists a major triple $E'=(E_1', E_2', E_3')$ of $(E, \beta)$ such that
\begin{equation}\label{VVgrtest}
|\int  \left(\sum |T_j(f_j,g_j)|^2 \right)^{\frac{1}{2}} (x)h(x) dx | \le K |E|^{\beta}
\end{equation}
for all sequences of functions  $ f = \{ f_{j}\} \in X_{\infty}(E'_1), g = \{ g_{j}\} \in X_{2}(E'_2), $ and $h \in X(E'_3)$---
\end{quotation}
then we have,
$$\| T(f,g) \|_{ L^r(l^2) }=\left\| \|\left\{ T_j(f_j,g_j) \right\}\|_{l^{2}} \right\|_r \lesssim K \left\| f  \right\|_{L^p(l^\infty)}\left\|  g \right\|_{L^q(l^2)},$$
where $(\frac{1}{p}, \frac{1}{q}, \frac{1}{r'} )= \alpha$ and $\frac{1}{r'}+ \frac{1}{r}=1.$
\end{lemma}

At this point we say a few words about how the previous lemmas will be used. The above vector valued estimate \eqref{VVgrtest} can be further reduced to showing estimates for single operators. More precisely, we will rely on the simple H\"older estimate
$$|\lb \left( \sum |T_j(f_j,g_j)|^2 \right)^{\frac{1}{2}}, h \rb | \lesssim  \left( \sum_j\int |T_j(f_j,g_j)|^2 1_{E'_3 } \right)^{\frac{1}{2}} |E_3|^{\frac{1}{2}}$$
which separates the contributions for individual $j$ terms. Thus, \eqref{VVgrtest} will follow if we can prove that, uniformly in $j$
\begin{equation}\label{Tiest}
|| T_j(f_j 1_{E'_1},g_j 1_{E'_2})1_{E'_3}||_2 \le K K' \|f_j\|_{\infty} \|g_j\|_2
\end{equation}
where
\begin{equation}
K'=K'_\beta=|E_1|^{\beta_1}|E_2|^{\beta_2- \frac{1}{2}} |E_3|^{\beta_3 -\frac{1}{2}}.
\end{equation}
Define the bilinear operator by
\begin{equation}\label{Sj}
S^{E'}_j(f,g) =T_j(f 1_{E'_1},g 1_{E'_2})1_{E'_3}
\end{equation}

To show \eqref{Tiest}, using Theorem \ref{grtthm}, it is enough to show $S^{E'}_j$ is generalized restricted type $\gamma \in \pi_0$ with constant $KK'$ for $\gamma$ in a neighborhood of the point $(0, \frac{1}{2}, \frac{1}{2}) \in \pi_0.$ Similar considerations apply to \eqref{VVrtest}. We summarize our progress so far as follows

\begin{lemma}\label{BTinter2l1} Let $\alpha \in \pi_1$ and  let $T =\{T_j \}$ be a sequence of bilinear operators. Assume there is a neighborhood $V$ of $\alpha$ in $\pi_1$ such that for each $\beta\in V$ and each  $E= (E_1, E_2, E_3)$ there exists a major triple $E' =(E'_1, E'_2, E'_3)$ of $(E, \beta)$ such that the associated trilinear forms for  $S^{E'}_j$ are generalized restricted type $\gamma$ for $\gamma$ in some neighborhood of the point $(0, \frac{1}{2}, \frac{1}{2})$, with constant $KK'_\beta$. Then we have
$$\| T(f,g) \|_{ L^r(l^2) }=\left\| \|\left\{ T_j(f_j,g_j) \right\}\|_{l^{2}} \right\|_{L^r} \lesssim K \left\| f  \right\|_{L^p(l^{\infty})}\left\|  g \right\|_{L^q(l^2)} ,$$
where $(\frac{1}{p}, \frac{1}{q}, \frac{1}{r'} )= \alpha$ and $\frac{1}{r'}+ \frac{1}{r}=1.$
\end{lemma}

\begin{lemma}\label{BTinter2} Let $\alpha \in \pi_0$ with first two indices good and let $T =\{T_j \}$ be a sequence of bilinear operators. Assume there is a neighborhood $V$ of $\alpha$ in $\pi_0$ such that for each $\beta\in V$ and each  $E= (E_1, E_2, E_3)$ there exists a major triple $E' =(E'_1, E'_2, E'_3)$ of $(E, \beta)$ such that the associated trilinear forms for  $S^{E'}_j$ are generalized restricted type $\gamma$ for $\gamma$ in some neighborhood of the point $(0, \frac{1}{2}, \frac{1}{2})$, with constant $KK'_\beta$. Then we have
$$\| T(f,g) \|_{ L^r(l^2) }=\left\| \|\left\{ T_j(f_j,g_j) \right\}\|_{l^{2}} \right\|_{L^r} \lesssim K \left\| f  \right\|_{L^p(l^{\infty})}\left\|  g \right\|_{L^q(l^2)} ,$$
where $(\frac{1}{p}, \frac{1}{q}, \frac{1}{r'} )= \alpha$ and $\frac{1}{r'}+ \frac{1}{r}=1.$
\end{lemma}

\section{Generalized restricted type interpolation for operators taking values in $L^P$ spaces with mixed norms}\label{grtLPsec}

Next we further generalize the interpolation Lemmas \ref{vvrt}, \ref{vvgrt} by replacing inner spaces $l^2, l^\infty$ by  $\l^P$ spaces with mixed norms \cite{BP1}, and  interpolate between those inner spaces. Interpolation in inner spaces is done using complex interpolation \cite{BP1} while interpolation in the outer spaces ($L^p$) is done using generalized restricted type interpolation. We start by  recalling  $\l^P$ spaces with mixed norms from \cite{BP1}.

\begin{definition}\label{LP}
Let $(X_i, \mu_i)$, for $1 \le i \le n$ be some  totally $\sigma$-finite measure spaces and let $P= (p_1, \dots , p_n) $ be some $n$-tuple with $1\le p_i \le \infty$. Let  $(X, \mu)= \prod X_i$ be the product space and let $f= f(x_1, x_2, \dots, x_n)$ be a measurable function on $X$. Define $\| f\|_P$ by taking $L^{p_i}$ norms successively,
$$ \| f\|_P = \|  \dots \| \dots \| f \|_{L^{p_1}(X_1)} \dots \|_{L^{p_i}(X_i)} \dots \|_{L^{p_n}(X_n)}$$
We denote the space of those $f$ with finite $\| f\|_P$ by $\l^P(X)$.
\end{definition}

In this section we use $P, Q, R$ to denote $n$-tuples as in  Definition \ref{LP}.  We use the notation $\vec{\alpha} = (\frac{1}{P}, \frac{1}{Q}, \frac{1}{R})$ to denote $\vec{\alpha} \in \pi_1^n$ where $\vec{\alpha} =(\alpha^{[1]}, \alpha^{[2]}, \dots, \alpha^{[n]})$ and $\alpha^{[i]} = (\frac{1}{p_i}, \frac{1}{q_i}, \frac{1}{r_i})$. Note that we have nested sequence of sets $I_m = A_{1,m}\times A_{2,m} \times \dots \times A_{n,m}$, with $A_{i,m} \subset X_i$ and  $\mu_i(A_{i,m}) < \infty $ and $ \cup_m I_m =X.$ Following the notations from \cite{BP1}, we denote by $H(X)$ the set of all functions of the form $\sum_k c_k 1_{B_k}$
where sum is over finitely many $k$ and $\cup_k B_k \subset I_m$ for some $m$ and $B_k$ is of the form $B_k =B_{1,k} \times B_{2,k} \times \dots \times B_{n,k}$.

We record the following fact  \cite{BP1} about $L^P$. Let $1 \le P < \infty$, that is $1\le p_i <\infty$ for $1\le i \le n$. Then we have
\begin{equation}\label{dualLP}
 \| f \|_P =\sup_{g\in H(X)} \int fg d\mu
\end{equation}
where the supremum is taken over $g$ with $\|g\|_{P'} \le1$.

Now we consider functions $F(x,t)=F(x_1, x_2, \dots, x_n, t)$ defined in $X \times \R^m$. For $0 < p \le \infty$ we denote $$\| F\|_{L^p(\l^P)} = \|  \| F\|_{\l^P} (t) \|_{L^p}.$$

Now we define the analogue of $X_P(E)$ in this settings. For $0 < P \le \infty$, that is $1\le p_i <\infty$ for $1\le i \le n$ and for a set $E$ with finite measure we set
\begin{equation}\label{XLP}
X_{\l^P}(E) = \{ F \in L^{\infty}(\l^P): F \in H(X \times \R^m), \| F\|_{\l^P} (t) \le 1_E(t)  \}.
\end{equation}

Now for a bilinear operator $T(F,G)(x,t)$ we define the associated trilinear form by,
$$\Lambda(F, G, H) = \int T(F, G)(x, t) H(x, t) d\mu(x) dt$$
\begin{definition}\label{grtLP}
If $(\vec{\alpha}, \alpha) \in  \pi_1^n \times \pi_0$, we say that a trilinear form $\Lambda$ is of  generalized restricted type $(\vec{\alpha}, \alpha)$ with constant $K$ if for triple $E= (E_1,E_2,E_3)$ there exists a major triple $E'$ of $(E, \alpha)$ such that
\begin{equation}\label{grtjgkfjkgjfkgj}
|\Lambda(F_1, F_2, F_3)|  \lesssim K |E|^{\alpha}
\end{equation}
for all functions $ F_i \in X_{\l^{P_i}}(E'_i), i = 1,2,3,$ where $\vec{\alpha} = (\frac{1}{P_1}, \frac{1}{P_2}, \frac{1}{P_3})$ and $|E|^{\alpha} = |E_1|^{\alpha_1}|E_2|^{\alpha_2}|E_3|^{\alpha_3}$.

If  $\Lambda$ is of  generalized restricted type $(\vec{\alpha}, \alpha)$ for  $\alpha \in \pi_1$ we say $\Lambda$ is of restricted type $(\vec{\alpha}, \alpha)$.
\end{definition}

Note that Remark \ref{noexcinrt} still holds in this situation.
Also note that we for $ F_i \in X_{\l^{P_i}}(E'_i), i = 1,2,$ by using \eqref{dualLP} we have,
\begin{equation}\label{relation}
 \sup_{F_3 \in  X_{\l^{P_3}}(E'_3)}| \Lambda(F_1, F_2, F_3) | \sim \int \| T(F_1, F_2)\|_{\l^{P'_3}}(t) 1_{E'_3}(t) dt.
\end{equation}
Compare this with the generalized restricted type condition in Lemma \ref{vvgrt}. Note also that Definition \ref{grtLP} allows us to use the linearity of $T$ in inner space, which allow us to use complex interpolation methods to interpolate between inner spaces. We prove the following version of  Lemma \ref{convexity}. The proof is a combination of the proof of the interpolation theorem from  \cite{BP1}, and the proof of the convexity Lemma \ref{convexity}. Also note that usually in complex interpolation we use strong endpoints but in here we have weak endpoint estimates, getting strong end points is done by next two real interpolation theorems, Theorems \ref{rtthmLP} and \ref{grtthmLP}.

\begin{lemma}\label{convexity2}
If $\Lambda$ is generalized restricted type $(\vec{\alpha}, \alpha)$ and $(\vec{\beta}, \beta)$  with constant $K_1$ and $K_2$ respectively,  with the same  choice of index $j$ and the major subset $E_j'$ , then $\Lambda$ is generalized restricted type $( [\vec{\alpha}, \vec{\beta}]_{\theta}, [\alpha, \beta]_{\theta})$ with the constant $K_1^{1-\theta} K_2^\theta$, where $[\alpha, \beta]_{\theta} =(1-\theta)\alpha  +\theta \beta $ and $\theta \in (0, 1)$, $[\vec{\alpha}, \vec{\beta}]_{\theta}$ is defined similarly.
\end{lemma}

Proof of  Lemma \ref{convexity2}: Let $\vec{\alpha} = (\frac{1}{P_1}, \frac{1}{Q_1}, \frac{1}{R_1})$, $\vec{\beta} = (\frac{1}{P_2}, \frac{1}{Q_2}, \frac{1}{R_2})$, $[\vec{\alpha}, \vec{\beta}]_{\theta}=\vec{\gamma} = (\frac{1}{P}, \frac{1}{Q}, \frac{1}{R})$, and $P=(p_1, p_2, \dots , p_n)$. Let $E=(E_1, E_2, E_3)$ be a triple, and $E'=(E'_1, E_2', E_3')$ be a triple given by the generalized restricted type conditions. Let $F_1 \in X_{\l^P}(E_1'), F_2 \in X_{\l^Q}(E_2')$, and $F_3 \in X_{\l^R}(E_3')$.

Consider the strip $S = \{z: 0 < Re(z) < 1 \}$. We define the following tuples by,
\begin{equation}\label{indexz}
P(z) = (1-z)\frac{P}{P_1}+z \frac{P}{P_2}, \\
Q(z) = (1-z)\frac{Q}{Q_1}+z \frac{Q}{Q_2}, \\
R(z) = (1-z)\frac{R}{R_1}+z \frac{R}{R_2} .
\end{equation}

Notice that $P(\theta)= Q(\theta)= R(\theta)= (1, \dots ,1) \in \Z^n$. Next we define a family of functions $F^z$ for $z\in \bar{S}$. First define the ``projections" $F^z_i$ of $F$,  for $1 \le i \le n-1$,
\begin{equation}\label{Fzi}
F^z_i(x_{i+1}, \dots, x_n,t)= \left( \| \dots \| F(x_1, \dots , x_i, x_{i+1}, \dots , x_n, t) \|_{L^{p_1}(X_1)} \dots \|_{L^{p_i}(X_i)}  \right)^{p_{i+1}(z)-p_1(z)}
\end{equation}
where $P(z) = (p_1(z), \dots, p_n(z)).$
Now define,
\begin{equation}\label{Fz}
F^z(x,t)= |F(x,t)|^{p_1(z)} (\text{Sign } F) \prod_{i=1}^{m-1} F^z_i.
\end{equation}

Note that $F^z \in H(X \times \R^m)$, $ F^{\theta} =F$,  $F^{iy} \in X_{\l^{P_1}}(E'_1)$, and $ F^{1+iy} \in X_{\l^{P_2}}(E'_1)$. Similarly define $G^z$ and $H^z$. Now define $\Psi(z)= \Lambda(F^z, G^z, H^z)$ using linearity of $\Lambda$ and the fact that $F^z, G^z, H^z \in H(X \times \R^m)$, we get that $\Psi$ is analytic in the strip $S$ and continuous in $\bar{S}$.

Now using the generalized restricted type $(\vec{\alpha}, \alpha)$ condition for $z \in \bar{S}$ with $Re(s)=0$ we get,
$$|\Psi(z)| = |\Lambda(F^z, G^z, H^z)| \le K_1 |E|^{\alpha}.$$
Similarly for $z \in \bar{S}$ with $Re(s)=1$ we get,
$$|\Psi(z)| = |\Lambda(F^z, G^z, H^z)| \le K_2 |E|^{\beta}.$$
Now using Hadamard's three lines lemma we get,
$$|\Psi(z)| = |\Lambda(F^z, G^z, H^z)| \le K_1^{1-\theta}K_2^{\theta} |E|^{[\alpha, \beta]_\theta}.$$
when  $Re(z) =\theta$. Now taking $z = \theta$ completes the proof of  Lemma \ref{convexity2}.

Next we state the interpolation theorems in this setting, these two theorems are about fixed inner spaces and they are simple generalizations of  generalized restricted type interpolation theorems, Theorem \ref{grtthm} and \ref{rtthm}.

\begin{thm}\label{rtthmLP} Let $( \vec{\alpha}, \alpha) \in \pi_1^n \times \pi_1$ and let $T$ be a bilinear operator. If the associated trilinear form $\Lambda_T$ is  restricted type $ ( \vec{\alpha} ,\beta )\in \pi_1^n \times \pi_1$ for $\beta \in \pi_1$ in  a neighborhood of $\alpha$ with constant $K,$ then we have
$$\|T(F,G)\|_{L^r(\l^R)} \lesssim K \|F\|_{L^p(\l^P)} \|g\|_{L^q(\l^Q)};$$
here $(\frac{1}{p}, \frac{1}{q}, \frac{1}{r'} )= \alpha$ and $\frac{1}{r'}+ \frac{1}{r}=1$ and $\vec{\alpha} =(\frac{1}{P}, \frac{1}{Q}, \frac{1}{R'})$.
\end{thm}

\begin{thm}\label{grtthmLP} Let $\alpha \in \pi_0\setminus \pi_1$ with first two good indices, $\vec{\alpha} \in \pi_1^n$, and let $T$ be a bilinear operator  which is restricted type $( \vec{\alpha}, \alpha_0)$ for some $\alpha_0 \in \pi_1$. If the associated trilinear form $\Lambda_T$ is generalized  restricted type $(\vec{\alpha}, \beta)$ for $\beta \in \pi_0$ in  a neighborhood of $\alpha$ with constant $K,$ then we have
$$\|T(F,G)\|_{L^r(\l^R)} \lesssim K \|F\|_{L^p(\l^P)} \|g\|_{L^q(\l^Q)};$$
here $(\frac{1}{p}, \frac{1}{q}, \frac{1}{r'} )= \alpha$ and $\frac{1}{r'}+ \frac{1}{r}=1$ and $\vec{\alpha} =(\frac{1}{P}, \frac{1}{Q}, \frac{1}{R'})$.
\end{thm}

First we show how to combine Lemma \ref{convexity2} with above Theorems. Let $( \vec{\alpha}, \alpha), (\vec{\beta}, \beta),$ $ ( \vec{\gamma}, \gamma) \in \pi^n_1 \times \pi_1$ and $(\vec{\delta}, \delta)$ be a
 point in the convex hull of   $\{ ( \vec{\alpha}, \alpha), (\vec{\beta}, \beta), ( \vec{\gamma}, \gamma) \}$ with convex hull of  $\{ \alpha, \beta, \gamma \}$ containing a neighborhood of $\delta$. If $T$ is generalized restricted type $( \vec{\alpha}, \alpha), (\vec{\beta}, \beta),$ and $ ( \vec{\gamma}, \gamma)$ then using Lemma \ref{convexity2} we get that $T$ is generalized restricted type $(\vec{\delta}, \delta)$. Having one such point is not enough to apply the Theorem \ref{grtLP}. But if we further assume that $T$ is generalized restricted type $( \vec{\alpha}, \alpha')$ for $\alpha'$ in a neighborhood of $\alpha$ and similarly for the other two points $(\vec{\beta}, \beta)$ and $( \vec{\gamma}, \gamma)$, then by Lemma \ref{convexity2}, we get that $T$ is generalized restricted type of $(\vec{\delta}, \delta')$ for $\delta'$ in a neighborhood of $\delta$, so we can apply the above theorems now.


We prove Theorem \ref{grtthmLP} now. Let $v_1= (-2, 1, 1), v_2= (1, -2, 1), v_3= (1, 1, -2)$. Using  Lemma \ref{convexity2} we get that there exists $\mu>0$ such that $T$ is generalized restricted type $(\vec{\alpha}, \alpha^{[i]})$ for $i =1,2,3 $ where $\alpha^{[i]} =\alpha + \mu v_i$. Denote $\vec{\alpha}=(\frac{1}{P}, \frac{1}{Q}, \frac{1}{R'})$ and $\alpha=(\frac{1}{p}, \frac{1}{q}, \frac{1}{r'})$.

Let $F, G \in H(X \times \R^m)$, assume that $\|F\|_{\l^P}(t)$, $\|G\|_{\l^Q}(t)$ is supported in $(0, \infty)$ and decreasing. For $k_1, k_2 \in \Z$ let $$F_{k_1}(x,t) = F(x,t) 1_{[2^{k_1}, 2^{k_1+1}]}$$ and $$G_{k_2}(x,t) = G(x,t) 1_{[2^{k_2}, 2^{k_2+1}]}.$$  We have since $r\le 1$,
\begin{align*}
\| T(F,G)\|_{L^r(\l^R)}^r &= \int\|  \sum_{k_1, k_2 \in \Z} T(F_{k_1},G_{k_2}) \|_{\l^R}^r(t) dt \\
 &\le   \sum_{k_1, k_2 \in \Z}   \int \|T(F_{k_1},G_{k_2}) \|_{\l^R}^r(t) dt
\end{align*}
Now assume  $\|T(F_{k_1},G_{k_2}) \|_{\l^R}(t)$ is supported in $(0, \infty)$ and decreasing.
So we get
\begin{align*}
\| T(F,G)\|_{L^r(\l^R)}^r &\le \sum_{k_1, k_2, k_3 \in \Z} \int \|T(F_{k_1},G_{k_2}) \|_{\l^R}^r(t) 1_{[2^{k_3}, 2^{k_3+1} ]}(t)dt\\
&\lesssim  \sum_{k_1, k_2, k_3 \in \Z}  2^{k_3(1-r)} \left( \int \|T(F_{k_1},G_{k_2}) \|_{\l^R}(t) 1_{[2^{k_3}, 2^{k_3+1} ]}(t)dt \right)^r.
\end{align*}
Given triple
\begin{equation}\label{setE}
E = ([2^{k_1}, 2^{k_1+1}], [2^{k_2}, 2^{k_2+1}] , [0, 2^{k_1+1}]),
\end{equation}
let $E_3'$ be the major set given by the generalized restricted type assumption, then we have
$$\lesssim  \sum_{k_1, k_2, k_3 \in \Z}  2^{k_3(1-r)} \left( \int \|T(F_{k_1},G_{k_2}) \|_{\l^R}(t) 1_{E_3'}(t)dt \right)^r.$$
Now using \eqref{dualLP} we get $H_{k_3} \in X_{\l^{R'}}(E_3')$ such that,
$$ \int \|T(F_{k_1},G_{k_2}) \|_{\l^R}(t) 1_{E_3'}(t)dt \lesssim \int T(F_{k_1}, G_{k_2})(x,t) H_{k_3}(x,t) dxdt$$
Consider trilinear form given by
$$\Lambda(F_{k_1}, G_{k_2}, H_{k_3}) = \int T(F_{k_1}, G_{k_2})(x,t) H_{k_3}(x,t) dxdt$$
Lets assume $k_1, k_2 \le k_3$, other two cases follows similarly. Now using the fact that $\Lambda$ is generalized restricted type of $(\vec{\alpha} , \alpha^{[3]})$, we get
$$| \Lambda(F_{k_1}, G_{k_2}, H_{k_3}) | \lesssim  |E|^{\alpha^{[3]}} \| F \|_{\l^P}(2^{k_1}) \| G \|_{\l^Q}(2^{k_2})  = |E|^{\alpha} |E|^{\mu v_3}  \| F \|_{\l^P}(2^{k_1}) \| G \|_{\l^Q}(2^{k_2}) .$$
now this with \eqref{setE} we get,
\begin{align*}
 &\sum_{k_1, k_2 \le k_3}  2^{k_3(1-r)} \left( \int \|T(F_{k_1},G_{k_2}) \|_{\l^R}(t) 1_{E_3'}(t)dt \right)^r\\
&\lesssim   \sum_{k_1, k_2 \le k_3} \| F \|^r_{\l^P}(2^{k_1})\| G \|^r_{\l^Q}(2^{k_2}) 2^{\frac{k_1r}{p}} 2^{\frac{k_2r}{q}} |E|^{\mu r v_3 }\\
&= \sum_{\bar{k_1}, \bar{k_2} \le 0} 2^{\mu \bar{k}_1 r} 2^{\mu \bar{k}_2 r}  \sum_{k_3} \| F \|^r_{\l^P}(2^{k_3 +\bar{k}_1})\| G \|^r_{\l^Q}(2^{k_3+\bar{k}_2}) 2^{\frac{(k_3+k_1)r}{p}} 2^{\frac{(k_3+k_2)r}{q}},
\end{align*}
here $\bar{k}_i=k_i -k_3$ for $i=1,2$. Now using the fact that $ \left( \sum 2^k  \| F\|^p_{\l^P}(2^k) \right)^{1/p} \lesssim \|F\|_{L^p(\l^P)}$ and similarly for $G$ we get,
$$\lesssim \sum_{\bar{k_1}, \bar{k_2} \le 0} 2^{\mu \bar{k}_1 r} 2^{\mu \bar{k}_2 r} \| F \|^r_{L^p(\l^P)} \| G \|^r_{L^q(\l^Q)} \lesssim \| F \|^r_{L^p(\l^P)} \| G \|^r_{L^q(\l^Q)}.$$ This completes the proof of the Theorem \ref{grtthmLP}.

\section{Discretization}\label{disc}

In this section we introduce the discrete model sums for the operator $BP$. First we follow the notation and comments from \cite{MTT2} and introduce rank $1$ collections of tri-tiles, which  are used in the decomposition of the bilinear Hilbert transform  \cite{LT1}, \cite{MTT2}. We use rank $1$ collections in the decomposition of $BP$ in the first coordinate. We also introduce rank $0$ collections of tri-tiles, which are used in the decomposition of the paraproduct, and we use these in the decomposition of $BP$ in the second coordinate.

\begin{definition}\label{grid}
Let $n\ge 1$ and $\sigma \in \{ 1, 2, 3\}^n$. We define the shifted $n$-dyadic mesh
$D=D^n_\sigma$ to be the collection of cubes of the form
$$D^n_\alpha=\{ 2^j(k+(0,1)^n+(-1)^j\sigma): j \in \Z, k \in \Z^n\}.$$
We define a shifted dyadic cube to be any member of a shifted $n$-dyadic mesh.
\end{definition}

Note that for any cube $Q$ there exists a shifted dyadic cube $Q'$ such that $Q \subset \frac{7}{10} Q'$ and $|Q'| \sim |Q|$.

\begin{definition}\label{sparse}
 A subset $D'$ of a shifted $n$-dyadic grid $D$ is called sparse if for any two cubes $Q,Q' \in D $ with $Q \not =Q'$ we have $|Q| < |Q'|$ implies $10^9l(Q) < l(Q')$ and $|Q| = |Q'|$ implies $10^9Q \cap 10^9Q' = \emptyset$.
\end{definition}

Note that any subset of shifted $n-$dyadic grid, say when $n \le 3$, can be split into $O(1)$ sparse subsets.

\begin{definition}[Tri-tile]\label{tile}
Let $\sigma=(\sigma_1, \sigma_2, \sigma_3) \in \{0, 1, 3 \}^3$, and let $1\le i \le 3$. An i-tile with shift $\sigma_i$ is a rectangle $P =I_P \times \omega_P$ with area 1 and with $I_P \in D^1_0, \omega_P \in D^1_{\sigma_i}$. A tri-tile with shift $\sigma$ is a triple $\vec{P}= (P_1, P_2, P_3)$ such that $P_i$ is a i-tile with shift $\sigma_i$, and  $I_{P_i} = I_{\vec{P}}$ is independent of $i$. The frequency cube $Q_{\vec{P}}$ of a tri-tile is defined to be $\prod^3_{i=1} \omega_{P_i}.$
\end{definition}

\begin{definition}
A set $\vec{\mathbf{P}}$ of tri-tiles is called sparse if all tri-tiles in $\vec{\mathbf{P}}$ have the same shift and the set $\{ Q_{\vec{P}}: \vec{P} \in \vec{\mathbf{P}} \}$ is sparse.
\end{definition}

\begin{definition}\label{order}
Let $P$ and $P'$ be tiles. We write $P'<P$ if  $I_{P'} \subsetneq I_P$ and $3 \omega_P \subset 3 \omega_{P'}$, and $P' \le P$ if $P'<P$ or $P'=P$. We write $P' \lesssim P$ if $I_{P'} \subset I_P$ and $10^7 \omega_P \subset 10^7 \omega_{P'}$. We write $P' \lesssim' P$ if $P' \lesssim P$ and $P' \not \le P$.
\end{definition}

\begin{definition}[Rank $1$ collection of tri-tiles]\label{rank1}
A collection $\vec{\mathbf{P}}$ of tri-tiles is said to  have rank $1$ if one has the following properties for all $\vec{P},\vec{P'}\in \vec{\mathbf{P}}$:

\begin{itemize}
  \item If $\omega_{P_j}= \omega_{P_j'}$ for some $j\in\{1, 2, 3\}$, then $\omega_{P_i}= \omega_{P_i'}$ for all $i\in\{1, 2, 3\}$.
  \item If $|I_{\vec{P'}}| < |I_{\vec{P}}|$ and $P'_j \le P_j$ for some $j= 1, 2, 3,$ then $P'_i \lesssim' P_i$, for all $i\in\{1,2,3\}\setminus\{j\}.$
\end{itemize}
\end{definition}

Next we define rank $0$ collection of tri-tiles. With some minor modifications we can use these collection of tri-tiles to decompose $BP$ in the second coordinate.

\begin{definition}[Rank $0$ collection of tri-tiles]\label{rank0}
A collection $\vec{\mathbf{P}}$ of tri-tiles is said to have rank $0$ if one has the following properties for all $\vec{P}\in \vec{\mathbf{P}} $:

\begin{itemize}
  \item  There is at most one $i =1, 2, 3$ such that $ 0 \in \omega_{P_i}$, this $i$ is called the overlapping index of $\vec{P}$.

  \item   For non overlapping indexes $i'$ of $\vec{P}$ we have
 $$ \omega_{P_{i'}} = \left( \frac{1}{|I_{\vec{P}}|}, \frac{2}{|I_{\vec{P}}|} \right).$$
 \end{itemize}
 \end{definition}

Note that any rank $0$ collection $\vec{\mathbf{P}}$ of tri-tiles can be split into three sub collections
$$\vec{\mathbf{P}} = \bigcup _{i=1}^{3}\vec{\mathbf{P}}^{[i]} $$
with each tri-tile  $ \vec{P} \in \vec{\mathbf{P}}^{[i]}$ having overlapping index $i$, when it exists. We  call the collection $\vec{\mathbf{P}}^{[i]}$ an $i$-overlapping collection of rank $0$ tri-tiles. Note that we are allowing the tiles where neither component is overlapping  to be included in any of the collections.

\begin{definition}[Wave packet on a tile]\label{wavepacket}
Let $P$ be a tile. A wave packet on $P$ is a smooth function $\varphi_P$ which has Fourier support in $\frac{9}{10} \omega_P$ and obeys the estimates
$$ |\frac{d^\alpha}{d x^{\alpha}}[e^{-ic(\omega_Px)}\varphi_P(x)]| \lesssim_{M,\alpha} |I_P|^{-\frac12-\alpha} \tilde{\chi}_{I_P}^M(x) $$
for all $M>0$ and $\alpha\ge 0$.
\end{definition}

\begin{definition}[Product tiles]\label{Product tri-tile }

For pair of tri-tiles $\vec{B} =(B_1, B_2, B_3), \vec{P}=(P_1, P_2, P_3)$ define the product tri-tile by $s=(s_1, s_2, s_3)$ where the components $s_i=B_i \times P_i$ are called product tiles. We call  a smooth function $\varphi_{s_i}$  a wave packet associated with the product tile $s_i$ with spatial interval $R_s=I_s\times J_s$ if $\varphi_{s_i}$ has the Fourier support in $\frac{9}{10} \omega_{B_i} \times \frac{9}{10} \omega_{P_i}$ and
$$|\partial_x^\alpha\partial_y^\beta[e^{-ic(\omega_{B_i})x}e^{-ic(\omega_{P_i})y}\varphi(x,y)] |\lesssim_{\alpha,\beta,M} |I_s|^{-\alpha-1/2}|J_s|^{-\beta-1/2}  \tilde{\chi}_{I_{s}}^M(x) \tilde{\chi}_{J_{s}}^M(y).$$

\end{definition}

\begin{definition}[Collections of product tri-tiles]\label{producttilecollection}
A collection of product tri-tiles, denoted by  $\SSS$, is called sparse if
$$\mathbf{B}_{\SSS}=\{ \vec{B} : \vec{B} \times \vec{P} \in \SSS \text{ for some } \vec{P} \}$$
is a rank 1 sparse collection of tri-tiles, and
$$\mathbf{P}_{\SSS}=\{ \vec{P} : \vec{B} \times \vec{P} \in \SSS \text{ for some } \vec{B} \}$$
is a rank 0 collection of tri-tiles.

We call a sparse collection $\SSS^{[i]}$ of product tri-tiles  $i$-overlapping for some $i=1, 2, 3$, if $\mathbf{P}_{\SSS^{[i]}}$ is a rank 0 collection of $i$-overlapping tri-tiles.

We call a collection $\SSS^{[i]}_j$ of product tri-tiles $i$-overlapping with scale $j$, if it is $i$-overlapping as defined above and for all tri-tiles  $\vec{P} \in \mathbb{P}_{\SSS^{[i]}_j}$ we have $|I_{\vec{P}}| = 2^{-j}$. Note that in this case we have $\omega_{P_{i'}} = (2^j, 2^{j+1})$ for $i' \not =i$.
\end{definition}

Note that any sparse collection $\SSS$ of product tri-tiles can be split into three sub collections,
\begin{equation}\label{decompS}
 \SSS = \bigcup_{i=1}^3\SSS^{[i]} =\bigcup_{i=1}^3 \bigcup_{j \in \Z} \SSS^{[i]}_{j}.
\end{equation}
Now we can define the  model operators for  $BP$.
\begin{definition}\label{ModelT}
For a sparse collection of product tri-tiles $\SSS$ we define the associated model operator $T=T_{\SSS}$ by
\begin{equation} \label {eq: modelsum}
T(f,g)(x,y) = \sum_{s\in \SSS} \frac{\epsilon_s}{|R_s|^{1/2}} \lb f,\varphi_{s_1}\rb \lb g,\varphi_{s_2}\rb \varphi_{s_3},
\end{equation}
where $|\epsilon_s | \le 1,$ and $\varphi_{s_i}$ is a wave packets associated with  $s_i$.

\end{definition}

By standard arguments \cite{MPTT}, \cite{MTT2}, Theorem \ref{Main-Theorem} can be reduced to the following theorem about model operators $T$.

\begin{thm}\label{modelsum} Given a finite sparse collection $\SSS$ of  product tri-tiles,
 the model operator  $T$ defined by \eqref{eq: modelsum} maps $L^p\times L^q $ into $  L^r$ for $p,q>1, \frac{1}{p} +\frac{2}{q} <2, \frac{1}{q} +\frac{2}{p} <2$ with $\frac{1}{p}+\frac{1}{q}=\frac{1}{r}$, and the bound is independent of the collection $\SSS$.
\end{thm}

\begin{definition} We define the model operators associated to $\SSS^{[i]}$ by
\begin{equation} \label {eq: modelsumi}
T^{[i]}(f,g)(x,y) = \sum_{s\in \SSS^{[i]}} \frac{\epsilon_s}{|R_s|^{1/2}} \lb f,\varphi_{s_1} \rb \lb g,\varphi_{s_2}\rb \varphi_{s_3},
\end{equation}
where $|\epsilon_s | \le 1.$

\end{definition}
Note that by \eqref{decompS} we have
 $$T=T^{[1]}+T^{[2]}+T^{[3]}$$
 In the following theorem we prove  that $T^{[i]} , i=1,2,3,$ satisfy different ranges of $L^p$ bounds.

\begin{thm}\label{thmTi}
The model operators $T^{[i]}$ defined by \eqref{eq: modelsumi} map $L^p\times L^q $ into $  L^r$  for  $\frac{1}{p}+\frac{1}{q}=\frac{1}{r}$, $p ,q >1$ satisfying the additional restrictions

(a) $\frac{1}{p}+\frac{2}{q}<2 ;$ for $T^{[1]}$

(b) $ \frac{2}{p}+\frac{1}{q}<2;$ for $T^{[2]}$

(c) $r> 2/3$, for $T^{[3]}$

\end{thm}

We prove this theorem in Section  \ref{proofofmainthm}. Note that since $T=T^{[1]}+T^{[2]}+T^{[3]}$, Theorem \ref{Main-Theorem} follows immediately from Theorem \ref{thmTi}.

We prove the following simpler version of Theorem \ref{thmTi} in Section \ref{l2}, which captures the key ideas without getting into technicalities.

\begin{thm}\label{l2thm} The model operator  $T$ defined by \eqref{eq: modelsum} maps $L^p\times L^q $ into $  L^r$ for $2<p,q, r'<\infty $ with $\frac{1}{p}+\frac{1}{q}=\frac{1}{r}$. \end{thm}

Next we introduce the model sums associated with the operators $BP_j$.

\begin{definition}\label{defTjmodel}
We define the single scale  model operators $T_j=T_j^{[1]}$  associated with $\SSS^{[1]}_j$ by
\begin{equation}
T_j(f,g)(x,y) = \sum_{s\in \SSS^{[1]}_j} \frac{\epsilon_s}{|R_s|^{1/2}} \lb f,\varphi_{s_1}\rb \lb g,\varphi_{s_2}\rb \varphi_{s_3},
\end{equation}
\end{definition}
Again by standard arguments one can reduce  Theorem \ref{VVBHTPl2} to the following theorem about  model sums $T_j$.

\begin{thm}\label{VVBHTTj} When $2<p, q, r' < \infty$ we have
\begin{equation}
 \|\left\{ T_j(f_j,g_j) \right\}_j\|_{L^r(l^2)} \lesssim  \left\| f  \right\|_{L^p(l^{\infty})}\left\|  g \right\|_{L^q(l^2)}.
\end{equation}

\end{thm}
We prove this theorem in section \ref{VVBHTPsec}.

\section{Some results from time frequency analysis }

Even though our model operator $T$ is defined with collections of product tri-tiles $\SSS$, it is enough to consider single scale collections $\SSS^{[1]}_j$, as in Definition \ref{producttilecollection}. This is in fact one of the key observations in this paper, which allow us to get bounds for $BP$. By fixing $j$, we focus attention on tri-tiles which depend only on one parameter, and hence we preserve the existence of a relation of order, which is missing in the case of variable $j$. This will be critical to our approach.

 We start by defining trees in $\SSS^{[1]}_j$.

\begin{definition}[Tree]\label{tree}
For $1 \le i \le 3 $ and  $s_T \in \SSS^{[1]}_j$, we define an $i$-tree with top $s_T$  to be a collection of product tri-tiles $T \subset \SSS^{[1]}_j$ such that for all $s=\vec{B} \times \vec{P} \in T$ we have
$$ B_i \le  B_{T,i}  \text{    and   } \vec{P}=\vec{P_T},$$
where $s_T= \vec{B_T} \times \vec{P_T}.$ We denote by $R_T=R_{s_{T}}$. We say that $T$ is a tree if it is an $i$-tree for some $1\le i \le 3$.
\end{definition}

\begin{definition}\label{size}
For a collection $S \subset \SSS^{[1]}_j$ and a function $f$ we define the $i$-$size(S,f)$ for  $i =1,2,3$ by
$$\text{i-}size (S,f)= \sup_{T\subset S} \left( \frac{1}{|R_{T}|} \sum_{s\in T} |\lb f, \varphi_{s_i} \rb |^2 \right)^{\frac{1}{2}}$$
where the supremum is taken over  all the $i'$-trees  $T \subset S$ for $i' \not =i.$
\end{definition}
We will drop the $i$ dependence, whenever the value of $i$ does not matter.

The following lemma is a simpler variant of the size lemma from \cite{LL1}.  Unlike our tri-tiles in $\SSS_j$, the ones in \cite{LL1} have variable orientations, so the proof of the following lemma is less technical.

\begin{lemma}[Size Lemma]\label{sizelemma}
Given $S\subset \SSS^{[1]}_j$ and $f$,  we can decompose $S=S_{big} \cup S_{small}$ so that $size( S_{small}, f )<\frac{1}{2} size(S, f)$ and $S_{big} = \bigcup_{T\in \mathbb{T}_{big}} T$ with
$$ \sum_{T\in \mathbb{T}_{big}}|R_T| \lesssim \|f\|^2_2 size(S, f)^{-2}.$$
here the constant is independent of $j$.
\end{lemma}

We have the following classical bound on size; we use this to get better bounds on size by removing exceptional sets, using the strong maximal function.

\begin{prop}[Single tree estimate]\label{ste}
Let $r>1$.
Let $T\subset \SSS^{[1]}_j$ be a tree, then we have
$$size(T) \lesssim_r \sup_{s\in T}\inf_{x\in R_s} M_r(f)(x)$$
uniformly over $j$.
\end{prop}

Once we have the size Lemma \ref{sizelemma} and single tree estimate in Proposition \ref{ste}, we can -for reader's convenience- sketch a proof of the following theorem, in the local $L^2$ regime.
\begin{thm}\label{BHTTjthm}The operator $T_j$ from Definition \ref{defTjmodel}
 maps $L^p \times L^q$ into $L^r$ uniformly in $j$, when $\frac{1}{p}+\frac{1}{q}=\frac{1}{r}$ and $p,q> 1, r>\frac{2}{3}.$
\end{thm}

We note that, in spite of being a two dimensional object,  each $T_j$ is essentially the model sum for a bilinear Hilbert transform. One of the crucial estimates behind the proof of Theorem \ref{BHTTjthm} is the following lemma, which follows by repeated applications of  Lemma \ref{sizelemma}.

\begin{lemma}\label{finalest} Let $S \subset \SSS^{[1]}_j$ and  for $i=1, 2, 3$ define $ \sigma_i = i$-$size(S, f_i)$. Then we have
$$\Lambda_S(f_1, f_2, f_3) \lesssim \sum_{\substack{ 2^{-n_1} \le \sigma_1 \\  2^{-n_2} \le \sigma_2\\  2^{-n_3} \le \sigma_3}}  2^{-n_1 -n_2 -n_3}
\min_{i=1, 2, 3} \left\{ 2^{2n_i} \|f_i\|^2_2 \right\}
$$
\end{lemma}

Now we use Lemma \ref{finalest} to show $\Lambda_S$ is  restricted type $\alpha$ for $\alpha \in \pi_2$. The key observation here is that we can do this without removing exceptional sets, that is we can take the major triple $E'=E.$

Given a triple $E=(E_1, E_2, E_3)$ take
\begin{equation}\label{EE'}
E' =E
\end{equation}
and $f_i \in X(E_i)$. We have $\sigma_i \lesssim 1$ by Lemma \ref{ste}, so by Lemma \ref{finalest},
$$\Lambda_S(f_1, f_2, f_3) \lesssim \sum_{\substack{ 2^{-n_1} \lesssim 1 \\  2^{-n_2} \lesssim 1\\  2^{-n_3} \lesssim 1}}  2^{-n_1 -n_2 -n_3}
\min_{i=1, 2, 3} \left\{ 2^{2n_i} |E_i| \right\}.$$
Now for $a=(a_1,a_2,a_3) \in \pi_2$, we have
$$\Lambda_S(f_1, f_2, f_3) \lesssim \sum_{\substack{ 2^{-n_1} \lesssim 1 \\  2^{-n_2} \lesssim 1\\  2^{-n_3} \lesssim 1}}  2^{-n_1 -n_2 -n_3} \left( 2^{2n_i} |E_i| \right)^{a_1} \left( 2^{2n_i} |E_i| \right)^{a_2}\left( 2^{2n_i} |E_i| \right)^{a_3} .$$
So we get
\begin{equation}
\label{sgyutd347r843iotktoyi76iu0-}
\Lambda_S(f_1, f_2, f_3) \lesssim |E|^{a}.
\end{equation}
This shows that the trilinear form associated to $T_j$ is  restricted type $a$ for $a\in \pi_2$, so by  Theorem \ref{rtthm} we get the $L^p$ boundedness of $T_j$ in the local $L^2$ range.

\section{Proof of the local $L^2$ case in  Theorem \ref{Main-Theorem}}\label{l2}
 Recall that in Section \ref{disc} we have reduced Theorem  \ref{Main-Theorem} to the proof of boundedness of model sum operators. We prove Theorem \ref{l2thm} here. Since $T=T^{[1]}+T^{[2]}+T^{[3]}$, it is enough to show that each of the $T^{[i]}$ satisfies Theorem \ref{l2thm}. By duality it is enough to show that
$$ |\Lambda^{[i]}(f_1,f_2,f_3)| = | \lb T^{[i]}(f_1,f_2),f_3 \rb | \lesssim ||f_1||_p ||f_2||_q ||f_3||_r.$$
Now by the symmetry of the forms $\Lambda^{[i]}$ and the range of exponents $\pi_2$, it is enough to prove $L^p$ estimates for $T^{[1]}$.

Next note that we have the decomposition  $$T^{[1]} =\sum_{j \in \Z} T_j.$$

The power of this decomposition comes from the fact that it is a Littlewood-Paley decomposition  in the second variable, that is for $s=\vec{B} \times \vec{P} \in \SSS^{[1]}_j$ we have $\omega_{P_3}=(2^j, 2^{j+1})$. This observation gives
\begin{equation}\label{Tjeq1}
\Pi^l_j T^{[1]} =T_j
\end{equation}

Also $\omega_{P_2} =(2^j, 2^{j+1})$ hence
\begin{equation}\label{Tjeq2}
T_j(f,g) = T_j(f, \Pi^l_jg)
\end{equation}
The Littlewood-Paley theorem implies
$$\|g\|_{p} \sim \left\| \left( \sum_j |\Pi^l_jg|^2 \right)^{1/2} \right\|_{p}.$$
and
$$\|T(f,g)\|_{r} \sim \left\| \left( \sum_j |T_j(f,g)|^2 \right)^{1/2} \right\|_{r}.$$
It thus suffices to prove
\begin{equation} \label {l2est}
\left\| \left( \sum_j |T_j(f,g_j)|^2 \right)^{1/2} \right\|_{r} \lesssim \|f \|_p \left\| \left( \sum_j |g_j|^2 \right)^{1/2} \right\|_{q},
\end{equation} for arbitrary $g_j$, and then apply this result with  $g_j= \Pi_j^lg.$

This will follow from Theorem \ref{VVBHTTj}, using a constant sequence of functions for $f$.

\section{Proof of Theorem \ref{VVBHTTj}.}\label{VVBHTPsec}

We use Lemma \ref{BTinter2l1}. Let $\alpha=(\frac1p,\frac1q,\frac1{r'})$ be an arbitrary point in $\pi_2$. First define

$$E'_2 = E_2 \setminus \left\{ x:M_{(1+ \epsilon)}1_{E_1}(x) \gtrsim  \left(\frac{|E_1|}{ |E_2|} \right)^{\frac{1}{1+\epsilon}} \right\},$$
and  $E_j'=E_j$ for $j\not=2$. Note that $E'_2$ is a major subset of $E_2$ and $E'= (E_1, E_2', E_3)$ is a major triple of $(E, \alpha)$. According to Lemma \ref{BTinter2l1}, it suffices to prove  that the bilinear operator $S_j^{E'}$  is generalized restricted type $\beta$ with constant $KK'$ where
$K'=K'_\alpha= |E_1|^{\alpha_1}|E_2|^{\alpha_2- \frac{1}{2}} |E_3|^{\alpha_3 -\frac{1}{2}}$ for each $\beta$ in a small enough neighborhood $U\subset \pi_0$ of $(0, 1/2, 1/2).$

By the symmetry of both  $K'_\alpha$ and of the trilinear form for $S_j^{E'}$ (defined in  \eqref{Sj}) in the second and third components (recall that the first index is overlapping) we can assume $|E_2| \ge |E_3|$.
Thus  we have
$$K'= \left( \frac{|E_1|}{|E_2|} \right)^{\frac{1}{p}} \left( \frac{|E_2|}{|E_3|} \right)^{\frac{1}{2} -\frac{1}{r'}} \ge \left( \frac{|E_1|}{|E_2|} \right)^{\frac{1}{p}}.$$
Also by using Theorem \ref{BHTTjthm} we have
$$||S_j^{E'}(f,g)||_2 \lesssim ||f||_{\infty} ||g||_2.$$
So we can further  assume that $|E_1| \le |E_2|,$ otherwise we are done.

Now it is enough to show $S_j^{E'}$ is generalized restricted type $\beta $ with constant
$K\left(\frac{|E_1 |}{|E_2|} \right)^{\frac{1}{p}},$
for $\beta\in U$. We will show slightly more, namely that the exponent $\frac1p$ can be replaced with any $0<\delta<1$.

Let $F=(F_1, F_2, F_3)$ be a triple of subsets with finite measure. Take
 $$F_1'=F_1 \setminus \left\{ x: M_{(1+ \epsilon)}1_{F_2}(x) \gtrsim  \left(\frac{|F_2|}{ |F_1|} \right)^{\frac{1}{1+\epsilon}}  \text{or } M_{(1+ \epsilon)}1_{F_3}(x) \gtrsim  \left(\frac{|F_3|}{ |F_1|} \right)^{\frac{1}{1+\epsilon}}  \right\}
$$
Note that $F'=(F_1', F_2, F_3)$ is a major subset of $(F, \beta)$ for $\beta \in \pi_0$ sufficiently close to $(0,1/2, 1/2)$.
Let $f_i \in X(F_i')$ for $i= 1, 2, 3$.

Using standard arguments involving decomposing our collection into further sub collections based on the position relative to  the exceptional sets (see for example \cite{MTT2}), we can assume that we have the following upper bounds for sizes relative to $f_i$,

$$\sigma_1\lesssim \left(\frac{|E_1|}{ |E_2|} \right)^{\frac{1}{1+\epsilon}} ,$$

 $$\sigma_2\lesssim \left(\frac{|F_2|}{ |F_1|} \right)^{\frac{1}{1+\epsilon}},$$
  and
  $$\sigma_3\lesssim \left(\frac{|F_3|}{ |F_1|} \right)^{\frac{1}{1+\epsilon}}.$$

\begin{figure}[h]
\psset{unit=1.5cm}
\begin{pspicture}(-2,1.7321)(6,-5.1961)
\multido{\n=0+.5,\nn=4+-.5,\nnn=0+-.8660}{4}{\psline(\n,\nnn)(\nn,\nnn)}
\multido{\n=0+1,\nn=2+.5,\nnn=-3.4641+.8660}{4}{\psline(\n,0)(\nn,\nnn)}
\multido{\n=4+-1,\nn=2+-.5,\nnn=-3.4641+.8660}{4}{\psline(\n,0)(\nn,\nnn)}
\multido{\n=4.2+-.5,\nn=4.7+-.5,\nnn=0+-.8660}{4}{\psline{<-}(\n,\nnn)(\nn,\nnn)}
\pscircle[linestyle=solid](2,-1.732){.07}

\psellipse[linestyle=solid](2,-1)(.19,.06)

\rput[l](2.1,-1.632){$U$}
\rput[l](2.2,-1){$V$}
\rput[l](1.97,-1.9){$\cdot \  b$}


\multido{\n=-.35+1.00,\nn=-.1+1.0}{4}{\psline{->}(\n,.6062)(\nn,.1732)}
\multido{\n=.15+.50,\nn=.4+.5,\nnn=-1.4722+-.8660,\nnnn=-1.0392+-.8660}{4}{\psline{->}(\n,\nnn)(\nn,\nnnn)}
\rput[l](5,0){$\alpha_3=1$}
\rput[l](4.5,-.8660){$\alpha_3=\frac{1}{2}$}
\rput[l](4,-1.732){$\alpha_3=0$}
\rput[l](3.5,-2.5981){$\alpha_3=-\frac{1}{2}$}
\rput[Br](0,-1.732){$\alpha_2=-\frac{1}{2}$}
\rput[Br](.5,-2.5981){$\alpha_2=0$}
\rput[Br](1,-3.4641){$\alpha_2=\frac{1}{2}$}
\rput[Br](1.5,-4.3301){$\alpha_2=1$}
\rput[b](-.5,.8660){$\alpha_1=1$}
\rput[b](.5,.8660){$\alpha_1=\frac{1}{2}$}
\rput[b](1.5,.8660){$\alpha_1=0$}
\rput[b](2.5,.8660){$\alpha_1=-\frac{1}{2}$}
\end{pspicture}\caption{$U$, $V$ and $b$.}\label{Loc2proof}
\end{figure}
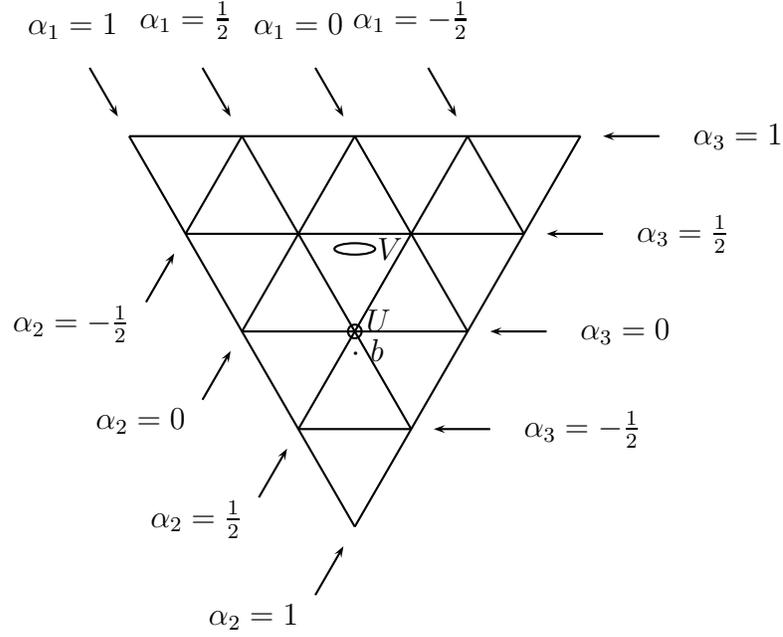

Now for $a=(a_1, a_2, a_3) \in \pi_2$, by Lemma \ref{finalest} we get
$$| \lb S_j^{E'}(f_1,f_2),f_3\rb | \lesssim \sum_{2^{-n_1} \le \sigma_1} \sum_{2^{-n_2} \le \sigma_2} \sum_{2^{-n_3} \le \sigma_3} 2^{-n_1 -n_2 -n_3} (2^{2n_1} |F_1|)^{a_1} $$
$$\cdot (2^{2n_2} |F_3|)^{a_3}(2^{2n_3} |F_3|)^{a_3} $$
\begin{equation}\label{sigmas}
\lesssim |F_1|^{a_1} |F_2|^{a_2}|F_3|^{a_3} \sigma_1^{1-2a_1} \sigma_2^{1-2a_2} \sigma_3^{1-2a_3}
\end{equation}
Choosing a triple $a\in \pi_2$ close to the point $(0,1/2, 1/2)$ and some $\epsilon $ sufficiently small we get
\begin{equation}\label{Sjl21}
| \lb S_j^{E'}(f_1,f_2),f_3\rb | \lesssim  \left(\frac{|E_1|}{|E_2|}\right)^{\delta(a)}|F_1|^{b_1} |F_2|^{b_2} |F_3|^{b_3},
\end{equation}
where $1>\delta(a)$ can be made as close to 1 as desired and $b=b(a)=(b_1, b_2, b_3) \in \pi_0$ can be made as close to the point $(0, 1/2, 1/2)$ as desired, with the additional restrictions $b_1<0, b_2 >1/2, b_3 >1/2$.

Using \eqref{sgyutd347r843iotktoyi76iu0-} we get that for each $d \in \pi_2$
\begin{equation}\label{Sjl22}
| \lb S_j^{E'}(f_1,f_2),f_3 \rb | \lesssim |F'|^{d}.
\end{equation}

The final part of the argument is rather delicate. Take  a small enough open set $V$ in $\pi_2$, far enough from the point $(0, 1/2, 1/2)$, and choose $a=(a_1, a_2, a_3) \in \pi_2$ as above, sufficiently close to $(0, 1/2, 1/2)$ such that the following holds: there is a neighborhood $U$ of $(0, 1/2, 1/2)$ which is contained inside the convex hull of $V$ and $b=b(a)$, such that each $\beta\in U$ can be represented as
$$\beta=\lambda b+(1-\lambda)c$$
with $c\in V$ and $\lambda\delta(a)>\frac1p$. See Figure \ref{Loc2proof}.

Let now $\beta=\lambda b+(1-\lambda)c\in U$.
Taking the geometric mean of $\eqref{Sjl21}$ and $\eqref{Sjl22}$ with $d=c$ we get
$$| \lb S_j^{E'}(f_1,f_2),f_3 \rb | \lesssim (|F'|^{c})^{1-\lambda}\left(\left(\frac{|E_1|}{|E_2|}\right)^{\delta(a)}|F|^b\right)^{\lambda}\lesssim \left(\frac{|E_1|}{|E_2|}\right)^{\frac1p}|F|^\beta,$$
completing the proof of Theorem \ref{VVBHTTj}.

\section{Proof of the main theorem: the full range}\label{proofofmainthm}
\begin{figure}[ht]
\psset{unit=.75cm}
\begin{pspicture}(-2,2.7321)(6,-5.1961)

\pspolygon[linestyle=dashed,fillcolor=gray,fillstyle=solid](.5,-.8660)(1,0)(1.5,-.8660)
\pspolygon[linestyle=dashed,fillcolor=gray,fillstyle=solid](2.5, -.866)(2,0)(3.5,-.866)
\pspolygon[linestyle=dashed,fillcolor=gray,fillstyle=solid](2.5, -2.658)(2,-1.732)(1,-1.732)

\multido{\n=0+.5,\nn=4+-.5,\nnn=0+-.8660}{4}{\psline(\n,\nnn)(\nn,\nnn)}
\multido{\n=0+1,\nn=2+.5,\nnn=-3.4641+.8660}{4}{\psline(\n,0)(\nn,\nnn)}
\multido{\n=4+-1,\nn=2+-.5,\nnn=-3.4641+.8660}{4}{\psline(\n,0)(\nn,\nnn)}
\multido{\n=4.2+-.5,\nn=4.7+-.5,\nnn=0+-.8660}{4}{\psline{<-}(\n,\nnn)(\nn,\nnn)}

\pspolygon[linestyle=dashed](.5,-.8660)(1,0)(1.5,-.8660)
\pspolygon[linestyle=dashed](2.5, -.866)(2,0)(3.5,-.866)
\pspolygon[linestyle=dashed](2.5, -2.658)(2,-1.732)(1,-1.732)

\multido{\n=-.35+1.00,\nn=-.1+1.0}{4}{\psline{->}(\n,.6062)(\nn,.1732)}
\multido{\n=.15+.50,\nn=.4+.5,\nnn=-1.4722+-.8660,\nnnn=-1.0392+-.8660}{4}{\psline{->}(\n,\nnn)(\nn,\nnnn)}
\rput[l](5,0){$\alpha_3=1$}
\rput[l](4.5,-.8660){$\alpha_3=\frac{1}{2}$}
\rput[l](4,-1.732){$\alpha_3=0$}
\rput[l](3.5,-2.5981){$\alpha_3=-\frac{1}{2}$}
\rput[Br](0,-1.732){$\alpha_2=-\frac{1}{2}$}
\rput[Br](.5,-2.5981){$\alpha_2=0$}
\rput[Br](1,-3.4641){$\alpha_2=\frac{1}{2}$}
\rput[Br](1.5,-4.3301){$\alpha_2=1$}
\rput[r]{-60}(-.5,.8660){$\alpha_1=1$}
\rput[r]{-60}(.5,.8660){$\alpha_1=\frac{1}{2}$}
\rput[r]{-60}(1.5,.8660){$\alpha_1=0$}
\rput[r]{-60}(2.5,.8660){$\alpha_1=-\frac{1}{2}$}
\end{pspicture}
\hspace{\stretch{1}}
\begin{pspicture}(-2,2.7321)(6,-5.1961)
\pspolygon[linestyle=dashed,fillcolor=gray,fillstyle=solid](1,-1.732)(1,0)(1.5,-.8660)
\pspolygon[linestyle=dashed,fillcolor=gray,fillstyle=solid](2.5, -.866)(3,0)(3,-1.732)
\pspolygon[linestyle=dashed,fillcolor=gray,fillstyle=solid](2.5, -2.6)(2,-1.732)(1.5,-2.6)

\multido{\n=0+.5,\nn=4+-.5,\nnn=0+-.8660}{4}{\psline(\n,\nnn)(\nn,\nnn)}
\multido{\n=0+1,\nn=2+.5,\nnn=-3.4641+.8660}{4}{\psline(\n,0)(\nn,\nnn)}
\multido{\n=4+-1,\nn=2+-.5,\nnn=-3.4641+.8660}{4}{\psline(\n,0)(\nn,\nnn)}
\multido{\n=4.2+-.5,\nn=4.7+-.5,\nnn=0+-.8660}{4}{\psline{<-}(\n,\nnn)(\nn,\nnn)}

\pspolygon[linestyle=dashed](1,-1.732)(1,0)(1.5,-.8660)
\pspolygon[linestyle=dashed](2.5, -.866)(3,0)(3,-1.732)
\pspolygon[linestyle=dashed](2.5, -2.6)(2,-1.732)(1.5,-2.6)

\multido{\n=-.35+1.00,\nn=-.1+1.0}{4}{\psline{->}(\n,.6062)(\nn,.1732)}
\multido{\n=.15+.50,\nn=.4+.5,\nnn=-1.4722+-.8660,\nnnn=-1.0392+-.8660}{4}{\psline{->}(\n,\nnn)(\nn,\nnnn)}
\rput[l](5,0){$\alpha_3=1$}
\rput[l](4.5,-.8660){$\alpha_3=\frac{1}{2}$}
\rput[l](4,-1.732){$\alpha_3=0$}
\rput[l](3.5,-2.5981){$\alpha_3=-\frac{1}{2}$}
\rput[Br](0,-1.732){$\alpha_2=-\frac{1}{2}$}
\rput[Br](.5,-2.5981){$\alpha_2=0$}
\rput[Br](1,-3.4641){$\alpha_2=\frac{1}{2}$}
\rput[Br](1.5,-4.3301){$\alpha_2=1$}
\rput[r]{-60}(-.5,.8660){$\alpha_1=1$}
\rput[r]{-60}(.5,.8660){$\alpha_1=\frac{1}{2}$}
\rput[r]{-60}(1.5,.8660){$\alpha_1=0$}
\rput[r]{-60}(2.5,.8660){$\alpha_1=-\frac{1}{2}$}
\end{pspicture}
\caption{Range of exponents for $\Lambda^{[2]}$ and $\Lambda^{[3]}$.}\label{exponentsT2T3}
\end{figure}
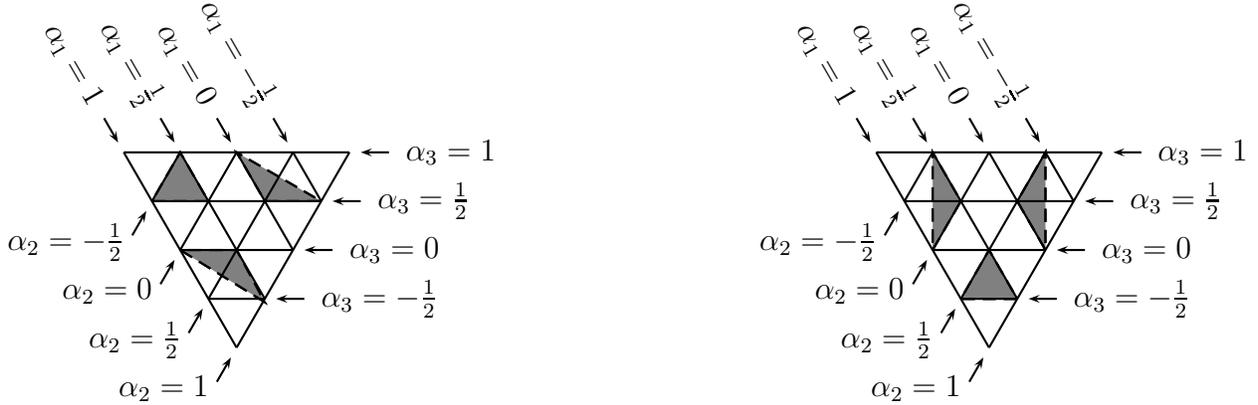

Consider the associated trilinear  form $\Lambda^{[i]}$ for $T^{[i]}$ given by
$$\Lambda^{[i]}(f_1,f_2,f_3) = \lb T^{[i]}(f_1,f_2),f_3 \rb .$$
We prove Theorem \ref{thmTi} by showing that $\Lambda^{[i]}$ is generalized restricted type $\alpha=(\alpha_1, \alpha_2, \alpha_3)$ for $ \alpha $ in the shaded triangles as shown in  Figures \ref{exponentspropT1} and \ref{exponentsT2T3}. Note that we have proved $\Lambda^{[i]}$ is restricted type $\alpha_0$ for $\alpha_0 \in \pi_2$ in Section \ref{l2}.

We prove the following Proposition to get Theorem \ref{thmTi}.
\begin{prop}\label{firstgrt}
Let  $\Lambda$ be the trilinear form with overlapping component 1 defined  by
$$\Lambda(f_1,f_2,f_3)= \sum_{s\in \SSS^{[1]}} \frac{\epsilon_s}{|R_s|^{1/2}} \lb f_1,\varphi_{s_1}\rb \lb f_2,\varphi_{s_2}\rb \lb f_3,\varphi_{s_3} \rb.$$
Then $\Lambda$ is generalized restricted type $(\alpha_1, \alpha_2, \alpha_3)$

(a) for  triples in the open triangle $A$ with vertices    $(-1/2,1,1/2)$, $ (-1/2,1/2,1)$ and $(0, 1/2, 1/2)$;

(b) for  triples  in the open triangle $B$ with vertices  $(1, -1/2, 1/2)$, $(1/2,0, 1/2)$ and $(0, 0, 1).$

\end{prop}

The triangles $A$ and $B$ are shown in the right and left shaded triangles respectively in the left picture in Figure \ref{exponentspropT1}.

 We first show how to get Theorem \ref{thmTi} using Proposition \ref{firstgrt}.  First note that in Proposition \ref{firstgrt} we have a trilinear form with overlapping first component and positive third component. Now for a $T^{[i]}$ and triple $\alpha$ in one of the shaded triangles corresponding to range of exponents for $T^{[i]}$ given in the Figures \ref{exponentspropT1} and \ref{exponentsT2T3}, consider the associated trilinear  form $\Lambda^{[i]}$ for $T^{[i]}$ given by,
$$\Lambda^{[i]}(f_1,f_2,f_3) = \lb  T^{[i]}(f_1,f_2),f_3 \rb.$$
Next  for a permutation $\sigma \in S_3$, define the trilinear form $\Lambda^{[i]}_{\sigma}$ by
 $$\Lambda^{[i]}_{\sigma}(f_1,f_2,f_3) = \Lambda^{[i]}(f_{\sigma(1)},f_{\sigma(2)},f_{\sigma(3)}).$$

One can choose a permutation $\sigma$ so that the $\Lambda^{[i]}_{\sigma}$ is having overlapping first component  and positive exponent in third component. Now use proposition \ref{firstgrt} to get the associated trilinear form $\Lambda^{[i]}$ is generalized restricted type $(\alpha_1, \alpha_2, \alpha_3)$. Finally note that  these give the ranges of exponents for $T^{[i]}$ as  in Theorem \ref{thmTi}.

\begin{figure}[ht]
\psset{unit=.75cm}

\begin{pspicture}(-2,2.7321)(6,-5.1961)
\pspolygon[linestyle=dashed,fillcolor=gray,fillstyle=solid](.5,-.8660)(2,0)(1.5,-.8660)
\pspolygon[linestyle=dashed,fillcolor=gray,fillstyle=solid](2.5, -.866)(3,0)(3.5,-.866)

\multido{\n=0+.5,\nn=4+-.5,\nnn=0+-.8660}{4}{\psline(\n,\nnn)(\nn,\nnn)}
\multido{\n=0+1,\nn=2+.5,\nnn=-3.4641+.8660}{4}{\psline(\n,0)(\nn,\nnn)}
\multido{\n=4+-1,\nn=2+-.5,\nnn=-3.4641+.8660}{4}{\psline(\n,0)(\nn,\nnn)}
\multido{\n=4.2+-.5,\nn=4.7+-.5,\nnn=0+-.8660}{4}{\psline{<-}(\n,\nnn)(\nn,\nnn)}
\pspolygon[linestyle=dashed](.5,-.8660)(2,0)(1.5,-.8660)
\pspolygon[linestyle=dashed](2.5, -.866)(3,0)(3.5,-.866)

\multido{\n=-.35+1.00,\nn=-.1+1.0}{4}{\psline{->}(\n,.6062)(\nn,.1732)}
\multido{\n=.15+.50,\nn=.4+.5,\nnn=-1.4722+-.8660,\nnnn=-1.0392+-.8660}{4}{\psline{->}(\n,\nnn)(\nn,\nnnn)}
\rput[l](5,0){$\alpha_3=1$}
\rput[l](4.5,-.8660){$\alpha_3=\frac{1}{2}$}
\rput[l](4,-1.732){$\alpha_3=0$}
\rput[l](3.5,-2.5981){$\alpha_3=-\frac{1}{2}$}
\rput[Br](0,-1.732){$\alpha_2=-\frac{1}{2}$}
\rput[Br](.5,-2.5981){$\alpha_2=0$}
\rput[Br](1,-3.4641){$\alpha_2=\frac{1}{2}$}
\rput[Br](1.5,-4.3301){$\alpha_2=1$}
\rput[r]{-60}(-.5,.8660){$\alpha_1=1$}
\rput[r]{-60}(.5,.8660){$\alpha_1=\frac{1}{2}$}
\rput[r]{-60}(1.5,.8660){$\alpha_1=0$}
\rput[r]{-60}(2.5,.8660){$\alpha_1=-\frac{1}{2}$}
\end{pspicture}
\hspace{\stretch{1}}
\begin{pspicture}(-2,1.7321)(6,-5.1961)

\pspolygon[linestyle=dashed,fillcolor=gray,fillstyle=solid](.5,-.8660)(2,0)(1.5,-.8660)
\pspolygon[linestyle=dashed,fillcolor=gray,fillstyle=solid](2.5, -.866)(3,0)(3.5,-.866)
\pspolygon[linestyle=dashed,fillcolor=gray,fillstyle=solid](1.5, -2.658)(2,-1.732)(3,-1.732)

\multido{\n=0+.5,\nn=4+-.5,\nnn=0+-.8660}{4}{\psline(\n,\nnn)(\nn,\nnn)}
\multido{\n=0+1,\nn=2+.5,\nnn=-3.4641+.8660}{4}{\psline(\n,0)(\nn,\nnn)}
\multido{\n=4+-1,\nn=2+-.5,\nnn=-3.4641+.8660}{4}{\psline(\n,0)(\nn,\nnn)}
\multido{\n=4.2+-.5,\nn=4.7+-.5,\nnn=0+-.8660}{4}{\psline{<-}(\n,\nnn)(\nn,\nnn)}

\pspolygon[linestyle=dashed](.5,-.8660)(2,0)(1.5,-.8660)
\pspolygon[linestyle=dashed](2.5, -.866)(3,0)(3.5,-.866)
\pspolygon[linestyle=dashed](1.5, -2.658)(2,-1.732)(3,-1.732)

\multido{\n=-.35+1.00,\nn=-.1+1.0}{4}{\psline{->}(\n,.6062)(\nn,.1732)}
\multido{\n=.15+.50,\nn=.4+.5,\nnn=-1.4722+-.8660,\nnnn=-1.0392+-.8660}{4}{\psline{->}(\n,\nnn)(\nn,\nnnn)}
\rput[l](5,0){$\alpha_3=1$}
\rput[l](4.5,-.8660){$\alpha_3=\frac{1}{2}$}
\rput[l](4,-1.732){$\alpha_3=0$}
\rput[l](3.5,-2.5981){$\alpha_3=-\frac{1}{2}$}
\rput[Br](0,-1.732){$\alpha_2=-\frac{1}{2}$}
\rput[Br](.5,-2.5981){$\alpha_2=0$}
\rput[Br](1,-3.4641){$\alpha_2=\frac{1}{2}$}
\rput[Br](1.5,-4.3301){$\alpha_2=1$}
\rput[r]{-60}(-.5,.8660){$\alpha_1=1$}
\rput[r]{-60}(.5,.8660){$\alpha_1=\frac{1}{2}$}
\rput[r]{-60}(1.5,.8660){$\alpha_1=0$}
\rput[r]{-60}(2.5,.8660){$\alpha_1=-\frac{1}{2}$}
\end{pspicture}
\caption{Range of exponents in proposition \ref{firstgrt} and $\Lambda^{[1]}$.}\label{exponentspropT1}
\end{figure}
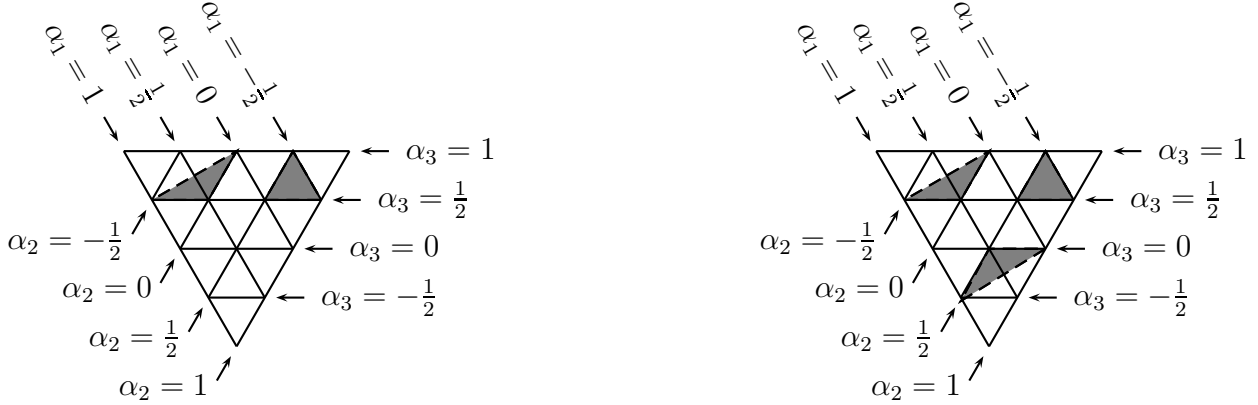

{\itshape{Proof of Proposition \ref{firstgrt}}:}  Proofs of the  cases (a) and (b) are similar so we do them simultaneously.

First note that by the symmetry of the second and third components in part (a) it is enough to consider the points from the open triangle $A'$  with vertices $(0, 1/2, 1/2)$, $( -1/2, 1, 1/2)$ and $(-1/4, 1/2, 3/4)$, note that we have  $\alpha_2 + 2\alpha_3<2 $ for  $ \alpha \in A'$.

Let $\alpha=(\alpha_1, \alpha_2, \alpha_3) \in A' \cup B$. Let $G=(G_1, G_2, G_3) \subset \R^2$ be a triple of measurable subsets with finite measure. We have to remove exceptional sets from the set corresponding to the negative exponent. We  use the strong maximal function $M$ as before, that is for $\alpha \in A'$ we remove an exceptional set from $G_1$ and when $ \alpha \in B$ from $G_2$.

In the first case, $\alpha \in A' $,  we take $$G_1'=G_1 \setminus \left\{ M_{1+ \epsilon}1_{G_2} \gtrsim  \left( \frac{|G_2|}{|G_1|} \right)^{\frac{1}{1+\epsilon}} \right\} $$ and take $G'=(G_1', G_2, G_3),$ and we prove the following estimate for $|f_i| \le 1_{G_i'} $:
\begin{equation} \label {eq: Trilinearform}
|\Lambda(f_1,f_2,f_3)|\lesssim |G|^{\alpha}= |G_1|^{\alpha_1} |G_2|^{\alpha_2} |G_3|^{\alpha_3}.
\end{equation}

In the second case, $\alpha \in B$, we remove an exceptional set from $G_2$ and define
$$G_2'=G_2 \setminus \left\{ M_{1+\epsilon}1_{G_1} \gtrsim  \left( \frac{|G_1|}{|G_2|} \right)^{\frac{1}{1+ \epsilon}} \right\}, $$
 and take $G'=(G_1', G_2, G_3)$, and we prove the same estimate \eqref{eq: Trilinearform} for  $|f_i| \le 1_{G_i'}$.

We undualize now: consider the  bilinear operator $T=T^{[1]}$ associated to  $\Lambda$ given by
$$T(f,g) = \sum_{s\in \SSS^{1}} \frac{\epsilon_s}{|R_s|^{1/2}} \lb f,\varphi_{s_1} \rb \lb g,\varphi_{s_2} \rb \varphi_{s_3};$$
so
 $$\Lambda(f_1,f_2,f_3) = \lb T(f_1,f_2),f_3 \rb.$$
Note that since $\alpha_3 \in (1/2,1)$ in both cases we have $ 1<\frac{1}{\alpha_3}<2.$ Let $r'=\frac{1}{\alpha_3}$ and let $r$ be the dual exponent of $r'$.  Notice to get \eqref{eq: Trilinearform} it is enough to show the  estimate
\begin{equation} \label {eq: estimateforT}
||T(f_1,f_2)||_{r} \lesssim |G_1|^{\alpha_1}|G_2|^{\alpha_2},
\end{equation}
To get \eqref{eq: estimateforT}  we view $T(f_1,f_2)$ as  a linear operator on the second component. For  $f_1$ fixed as above, consider, for the part (a), $\tilde{T}$ given by
\begin{equation} \label {eq: definitionofT1}
\tilde{T}(f_1,g) = \sum_{s\in \SSS^{1}} \frac{1}{|R_s|^{1/2}} \lb f_1,\varphi_{s_1} \rb \lb g1_{G_2},\varphi_{s_2}\rb \varphi_{s_3}.
\end{equation}
and for the part (b) consider,
\begin{equation} \label {eq: definitionofT1b}
\tilde{T}(f_1,g) = \sum_{s\in \SSS^{1}} \frac{1}{|R_s|^{1/2}} \lb f_1,\varphi_{s_1} \rb \lb g1_{G_2'},\varphi_{s_2}\rb \varphi_{s_3}.
\end{equation}
We will prove the following estimate for $\tilde{T}$:
$$||\tilde{T}(f_1,g)||_{r} \lesssim |G_1|^{\alpha_1}|G_2|^{\alpha_2 -1/{r}} ||g||_{r}.$$
Since $\alpha_2 -1/{r} =-\alpha_1$ we can write above estimate as
\begin{equation} \label {eq: estimateforT1}
||\tilde{T}(f_1,g)||_{r} \lesssim  \left(\frac{|G_1|}{|G_2|}\right)^{\alpha_1} ||g||_{r}.
\end{equation}
Now define linear ( $f_1$ is fixed) single scale operators  $T_j$ associated to the collections $\SSS^{[1]}_j$.  For part (a) define $T_j$ by
$$T_j(f_1,g)(x,y) = \sum_{s\in \SSS^{[1]}_j} \frac{\epsilon_s}{|R_s|^{1/2}} \lb f_1,\varphi_{s_1} \rb \lb g1_{G_2},\varphi_{s_2}\rb \varphi_{s_3},$$
and define $T_j$ for part (b) by
$$T_j(f_1,g)(x,y) = \sum_{s\in \SSS^{[1]}_j} \frac{\epsilon_s}{|R_s|^{1/2}} \lb f_1,\varphi_{s_1} \rb \lb g1_{G_2'},\varphi_{s_2}\rb \varphi_{s_3}.$$

Since we are working in the collection with lacunary second and third components, by using the classical (linear) Littlewood Paley theorem we get the estimates
$$||T(f_1,g)||_{r} \sim \left\| \left( \sum_j |T_j(f_1,g)|^2 \right)^{1/2} \right\|_{r},$$
$$||g||_{r} \sim \left\| \left( \sum_j |g_j|^2 \right)^{1/2} \right\|_{r},$$ where $g_j = \Pi^l_j g$.  Note also that we have $T_j(f_1, g) =T_j(f_1, g_j)$. So to get \eqref{eq: estimateforT1} it is enough to show that
\begin{equation} \label {eq: sqfunctionestimate}
\left\| \left( \sum_j |T_j(f_1,g_j)|^2 \right)^{1/2} \right\|_{r} \lesssim  \left(\frac{|G_1|}{|G_2|}\right)^{\alpha_1} \left\| \left( \sum_j |g_j|^2 \right)^{1/2} \right\|_{r}.
\end{equation}

We follow the proof of the  Lemma \ref{BTinter2} for the linear operators $T_j$ to reduce the vector valued estimates to estimates for single operators. Let $$K_0=\left(\frac{|G_1|}{|G_2|}\right)^{\alpha_1}.$$
We start with the following linear version of Lemma \ref{vvgrt}.

\begin{prop}\label{interpolationprop}

Suppose for any measurable sets $F_2,F_3 \subset \R^2$  with finite measure, we have $F_2' \subset F_2$ with $|F_2'| \ge \frac{1}{2} |F_2|,$ such that for any $ \left( \sum_j |g_j|^2 \right) \le 1_{F_2'}$ we have
\begin{equation} \label {eq: sqestimate1}
 \int \left( \sum_j \left|  T_j(f_1,g) \right|^2 \right)^{1/2} 1_{F_3} \lesssim  K_0|F_2|^{1/{\tilde{r}}}|F_3|^{1/{\tilde{r}'}}
\end{equation}
for $\tilde{r}$ in a small neighborhood of $r$. Then
$$ \left\| \left( \sum_j |T_j(f_1,g)|^2 \right)^{1/2} \right\|_{r} \lesssim  K_0 \left\| \left( \sum_j |g_j|^2 \right)^{1/2} \right\|_{r}.$$
\end{prop}
For $\tilde{r}$  consider $\tilde{\alpha_2}$ given by $ \tilde{\alpha_2} =1- \frac{1}{\tilde{r}'}-\alpha_1$. Note that when $\tilde{r}$ is close enough to $r$, the triple $(\alpha_1, \tilde{\alpha_2}, \frac{1}{\tilde{r}'})$ is in the same open triangle ($A'$ or $B$), as $\alpha$. For the rest of the argument we denote  $(\alpha_1, \tilde{\alpha_2}, \frac{1}{\tilde{r}'})$ simply by $\alpha=(\alpha_1, \alpha_2, \alpha_3)$ and $\tilde{r}$ by $r$.

Given $F_2$ and $F_3$ remove an exceptional set from $F_2$ and define
$$F_2'=F_2 \setminus  \left\{ M_{1+ \epsilon}1_{F_3} \ge \left( \frac{|F_3|}{|F_2|} \right)^{\frac{1}{1+\epsilon}} \right\} .$$
Note that since
$$\int \left( \sum_j \left|  T_j(f_1,g) \right|^2 \right)^{1/2} 1_{F_3}  \lesssim \left( \int \sum_j \left|  T_j(f_1,g) \right|^2  1_{F_3} \right)^{1/2} |F_3|^{1/2},$$
to get \eqref{eq: sqestimate1}  it is enough to show
\begin{equation} \label {eq: intL2estimate1}
\int |T_j(f_1,g)|^2 1_{F_3} \lesssim  K_0^2 |F_2|^{2/{r} -1}|F_3|^{2/{r'}-1}  ||g||^2_2.
\end{equation}
Note also that $-(2/{r} -1)= 2/{r'}-1=2\alpha_3-1$.

Now consider the bilinear operators $S_j=S_j^{F', G'}$ corresponding to the single scale collections $\SSS^{[1]}_j$. For the case (a) we define $S^j(f,g)$ by
\begin{equation} \label {eq: singlescaleoperator_a}
S^j(f,g) =  \sum_{s\in \SSS_j} \frac{1}{|R_s|^{1/2}} \lb f1_{G_1'},\varphi_{s_1} \rb \lb g1_{G_2 \cap F_2'},\varphi_{s_2}\rb \varphi_{s_3} 1_{F_3},
\end{equation}
and in  case (b) we define $S^j(f,g)$ by
\begin{equation} \label {eq: singlescaleoperator_b}
S^j(f,g) =  \sum_{s\in \SSS_j} \frac{1}{|R_s|^{1/2}} \lb f1_{G_1},\varphi_{s_1} \rb \lb g1_{G_2' \cap F_2'},\varphi_{s_2} \rb \varphi_{s_3} 1_{F_3}.
\end{equation}
Now to get \eqref{eq: intL2estimate1} it is enough to show the following $L^2$ estimate for the operators $S_j$:
\begin{equation} \label {eq: L2estimate3}
||S^j(f,g)||_2 \lesssim  K_0 K ||f||_{\infty}||g||_2,
\end{equation} with constant independent of $j,$
where $$K=  \left(\frac{|F_3|}{|F_2|}\right)^{\alpha_3-1/2}.$$

By generalized restricted type interpolation, to get \eqref{eq: L2estimate3} it is enough to show the associated trilinear form for $S^j(f,g)$ is generalized restricted type $t=(t_1,t_2,t_3)\in \pi_0$ in a neighborhood of  $(0,1/2,1/2)$, which  we prove in the following Proposition.

\begin{prop}\label{secondgrt}
Suppose for  given $E_1, E_2, E_3 \subset \R^2$ with finite measure  we have  $E_1'$ is a major set of $E_1$ such that for all  $|f| \le 1_{E'_1}, |g| \le 1_{E_2}$, and $|h| \le 1_{E_3}$,
\begin{equation} \label {eq: finalestimate}
|\lb S_j(f,g),h \rb | \lesssim K_0 K |E_1|^{t_1} |E_2|^{t_2}|E_3|^{t_3}
\end{equation}
for $t_1,t_2 \in B(1/2, \epsilon') $ with $ t_1+t_2+t_3=1,$ for sufficiently small $\epsilon' >0,$ and
with constant in \eqref{eq: finalestimate}  independent of $j$.
\end{prop}

Given $E=(E_1, E_2, E_3),$ we define $E_1'$ by
$$E_1'=E_1 \setminus  \left( \left\{ M_{1+\epsilon}1_{E_2} \ge \left( \frac{|E_2|}{|E_1|} \right)^{\frac{1}{1+\epsilon}}  \right\} \cup \left\{ M_{1+\epsilon}1_{E_3} \ge \left( \frac{|E_3|}{|E_1|} \right)^{\frac{1}{1+ \epsilon}} \right\} \right),$$
here $M$ denotes the strong maximal function as before; note also that the choice of $E'_1$ is independent of $j$ as needed.

For the rest of the proof we consider the two parts of Proposition \ref{firstgrt} separately. We first prove part (a). Using Lemma \ref{ste} and classical arguments involving decomposing $\SSS_j$ in to further sub collections based on the position relative  to the exceptional sets we can assume that we have following bounds on $\sigma_i$ for $i=1, 2,3$.
$$\sigma_1 \lesssim 1,$$
$$ \sigma_2 \lesssim \text{min} \left\{ \left( \frac{|G_2|}{|G_1|} \right)^{\frac{1}{1+\epsilon}}   , \left( \frac{|E_2|}{|E_1|} \right)^{\frac{1}{1+\epsilon}}  \right\},$$
and $$\sigma_3\lesssim  \text{min} \left\{ \left( \frac{|F_3|}{|F_2|} \right)^{\frac{1}{1+\epsilon}} , \left( \frac{|E_3|}{|E_1|} \right)^{\frac{1}{1+\epsilon}}  \right\}.$$
Now by using Proposition \ref{finalest} we get  for  $(\alpha, \beta, \gamma) \in \pi_2$,
$$|\lb S^j(f,g),h\rb | \lesssim |E_1|^{\alpha } |E_2|^{\beta} |E_3|^{\gamma}\sigma_1^{1-2\alpha}\sigma_2^{1-2\beta}  \sigma_3^{1-2\gamma}.$$
Now we have for $\nu ,\theta \in (0, \frac{1}{1+ \epsilon}),$
\begin{multline*}
| \lb S^j(f,g),h \rb | \lesssim |E_1|^{\alpha } |E_2|^{\beta} |E_3|^{\gamma}  \\ \cdot \left(\frac{|G_2|}{|G_1|}\right)^{\theta(1-2\beta)} \left(\frac{|E_2|}{|E_1|}\right)^{(1-\theta)(1-2\beta)}  \left(\frac{|F_3|}{|F_2|}\right)^{\nu(1-2\gamma)}\left(\frac{|E_3|}{|E_1|}\right)^{(1-\nu)(1-2\gamma)}.
\end{multline*}

Note since $\alpha \in A'$, we have $1- \alpha_3 \in (0,1/2)$, $1/2+\alpha_1\in (0,1/2)$, and  $2\alpha_3+\alpha_2 -\frac{3}{2} \in(0,1/2)$. Now choose  $\epsilon'$ small enough so that when $|\mu_2| ,|\mu_3| < \epsilon'$, we get $\gamma := 1-\alpha_3 -\mu_3 \in (0,1/2)$, $\beta:= 1/2+\alpha_1 -\mu_2 \in (0,1/2)$, and $\alpha := 1-\beta -\gamma =2\alpha_3+\alpha_2 -\frac{3}{2} +\mu_2 + \mu_3 \in (0,1/2).$

Finally note that $0 < -\alpha_1 <1-2\beta$ and $0<\alpha_3-\frac{1}{2}< 1-2\gamma$. So we can get $ \theta, \nu \in (0,1)$,  hence in  $ (0, \frac{1}{1+ \epsilon})$, so that  $\theta(1-2\beta) = -\alpha_1$ and $\nu(1-2\gamma) =\alpha_3-\frac{1}{2}$  With these choices of $\alpha, \beta, \gamma, \theta,$ and $ \nu$  we get
$$
|\lb S(f,g),h \rb | \lesssim \left(\frac{|G_1|}{|G_2|}\right)^{\alpha_1} \left(\frac{|F_3|}{|F_2|}\right)^{\alpha_3-1/2} |E_1|^{\mu_1} |E_2|^{1/2 +\mu_2}|E_3|^{1/2+\mu_3}
$$
here $\mu_1 +\mu_2+\mu_3=0$.  This completes the proof of part (a) of the proposition \ref{firstgrt}.

Next we  prove part (b). By similar argument as in part (a) we get, for $(\alpha, \beta, \gamma) \in \pi_2$
$$|\lb S^j(f,g),h \rb | \lesssim |E_1|^{\alpha } |E_2|^{\beta} |E_3|^{\gamma}\sigma_1^{1-2\alpha}\sigma_2^{1-2\beta}  \sigma_3^{1-2\gamma},$$
and following bounds on $\sigma_i$'s,
$$\sigma_1 \lesssim  \left( \frac{|G_1|}{|G_2|} \right)^{\frac{1}{1+\epsilon}} ,$$
$$ \sigma_2 \lesssim \text{min} \left\{ 1, \left( \frac{|E_2|}{|E_1|} \right)^{\frac{1}{1+\epsilon}} \right\},$$
and $$\sigma_3 \lesssim \text{min} \left\{ \left( \frac{|F_3|}{|F_2|} \right)^{\frac{1}{1+\epsilon}}, \left( \frac{|E_3|}{|E_1|} \right)^{\frac{1}{1+\epsilon}} \right\}.$$
For $ \nu ,\delta \in (0,\frac{1}{1+ \epsilon}),$ we have
\begin{multline*}
|\lb S^j(f,g),h \rb |  \lesssim |E_1|^{\alpha } |E_2|^{\beta} |E_3|^{\gamma} \\ \cdot  \left(\frac{|G_1|}{|G_2|}\right)^{(1-2\alpha)}  \left(\frac{|E_2|}{|E_1|}\right)^{\delta(1-2\beta)}   \left(\frac{|F_3|}{|F_2|}\right)^{\nu(1-2\gamma)}\left(\frac{|E_3|}{|E_1|}\right)^{(1-\nu)(1-2\gamma)}.
\end{multline*}

Now since $\alpha \in B$, we have that $1- \alpha_3 \in (0,1/2)$, $\frac{1}{2}-\frac{1}{2}\alpha_1  \in (0,1/2)$, and  $\alpha_3+ \frac{1}{2} \alpha_1 -\frac{1}{2} \in(0,1/2)$. Now for  $\epsilon'$ small enough and $|\mu_3| < \epsilon'$, we get $\gamma := 1- \alpha_3 -\mu_3 \in (0,1/2)$, $\alpha:= \frac{1}{2}-\frac{1}{2}\alpha_1  \in (0,1/2)$, and $\beta := 1-\gamma -\alpha =\alpha_3+ \frac{1}{2}\alpha_1 -\frac{1}{2}+\mu_3 \in (0,1/2).$

%

Finally note that for small enough $\epsilon'>0$ and $|\mu_2| <\epsilon',$ we get $0 < \alpha_3-\frac{1}{2} <1-2\gamma$ and $0<\frac{1}{2}-\beta +\mu_2< 1-2\beta$. Now we can get $  \nu, \delta \in (0,1) $ so that  $\nu(1-2\gamma) =\alpha_3-\frac{1}{2}$ and $\delta (1-2\beta)=\frac{1}{2}-\beta +\mu_2$.  With these choices of $\alpha, \beta, \gamma, \nu,$ and $ \delta$  we get
 $$
 |\lb S(f,g),h \rb | \lesssim \left(\frac{|G_1|}{|G_2|}\right)^{\alpha_1} \left(\frac{|F_3|}{|F_2|}\right)^{\alpha_3-1/2} |E_1|^{\mu_1} |E_2|^{1/2 +\mu_2}|E_3|^{1/2+\mu_3}.$$
This completes the proof of the proposition.

\section{Vector valued Bilinear Hilbert transform}\label{VVBHTsec}

In this section we prove the Theorems \ref{VVBHTPmess} and \ref{VVBHTPmess1}. One can follow the same line of arguments to prove the Theorem \ref{VVBHT}.

\begin{figure}[htbp]
\psset{unit=.75cm}

\begin{pspicture}(-2,2.7321)(4,-5.1961)

\pspolygon[linestyle=dashed,fillcolor=gray,fillstyle=solid](2,-2.0785)(2.6,-1.0392)(2.9,-1.5588)


\multido{\n=0+.5,\nn=4+-.5,\nnn=0+-.8660}{4}{\psline(\n,\nnn)(\nn,\nnn)}
\multido{\n=0+1,\nn=2+.5,\nnn=-3.4641+.8660}{4}{\psline(\n,0)(\nn,\nnn)}
\multido{\n=4+-1,\nn=2+-.5,\nnn=-3.4641+.8660}{4}{\psline(\n,0)(\nn,\nnn)}

\multido{\n=-.35+1.00,\nn=-.1+1.0}{4}{\psline{->}(\n,.6062)(\nn,.1732)}
\multido{\n=.15+.50,\nn=.4+.5,\nnn=-1.4722+-.8660,\nnnn=-1.0392+-.8660}{4}{\psline{->}(\n,\nnn)(\nn,\nnnn)}
\rput[Br](0,-1.732){$\alpha_2=-\frac{1}{2}$}
\rput[Br](.5,-2.5981){$\alpha_2=0$}
\rput[Br](1,-3.4641){$\alpha_2=\frac{1}{2}$}
\rput[Br](1.5,-4.3301){$\alpha_2=1$}
\rput[r]{-60}(-.5,.8660){$\alpha_1=1$}
\rput[r]{-60}(.5,.8660){$\alpha_1=\frac{1}{2}$}
\rput[r]{-60}(1.5,.8660){$\alpha_1=0$}
\rput[r]{-60}(2.5,.8660){$\alpha_1=-\frac{1}{2}$}
\end{pspicture}
\hspace{\stretch{.5}}
\begin{pspicture}(0,2.7321)(4,-5.1961)

\pspolygon[linestyle=dashed,fillcolor=gray,fillstyle=solid](1.2,-0.6928)(1.4,-1.0392)(2.6,-1.0392)(2.1,-0.1732)

\multido{\n=0+.5,\nn=4+-.5,\nnn=0+-.8660}{4}{\psline(\n,\nnn)(\nn,\nnn)}
\multido{\n=0+1,\nn=2+.5,\nnn=-3.4641+.8660}{4}{\psline(\n,0)(\nn,\nnn)}
\multido{\n=4+-1,\nn=2+-.5,\nnn=-3.4641+.8660}{4}{\psline(\n,0)(\nn,\nnn)}



\end{pspicture}
\hspace{\stretch{.5}}
\begin{pspicture}(0,2.7321)(6,-5.1961)

\pspolygon[linestyle=dashed,fillcolor=gray,fillstyle=solid](2,-2.0785)(2.9,-1.5588)(2.1,-0.1732)(1.2,-0.6928)


\multido{\n=0+.5,\nn=4+-.5,\nnn=0+-.8660}{4}{\psline(\n,\nnn)(\nn,\nnn)}
\multido{\n=0+1,\nn=2+.5,\nnn=-3.4641+.8660}{4}{\psline(\n,0)(\nn,\nnn)}
\multido{\n=4+-1,\nn=2+-.5,\nnn=-3.4641+.8660}{4}{\psline(\n,0)(\nn,\nnn)}
\multido{\n=4.2+-.5,\nn=4.7+-.5,\nnn=0+-.8660}{4}{\psline{<-}(\n,\nnn)(\nn,\nnn)}

%

\rput[l](5,0){$\alpha_3=1$}
\rput[l](4.5,-.8660){$\alpha_3=\frac{1}{2}$}
\rput[l](4,-1.732){$\alpha_3=0$}
\rput[l](3.5,-2.5981){$\alpha_3=-\frac{1}{2}$}
\end{pspicture}

\caption{Range of exponents  $A^{R}_1$, $A^{R}_2$, $A^{R}$ for $\frac{1}{R}=.6 >\frac{1}{2}$.}\label{ARR<2}
\end{figure}
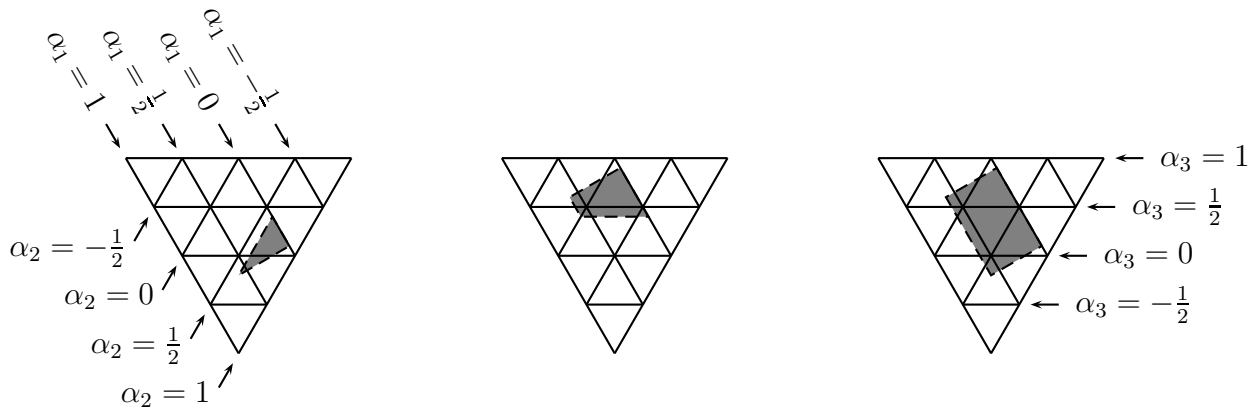

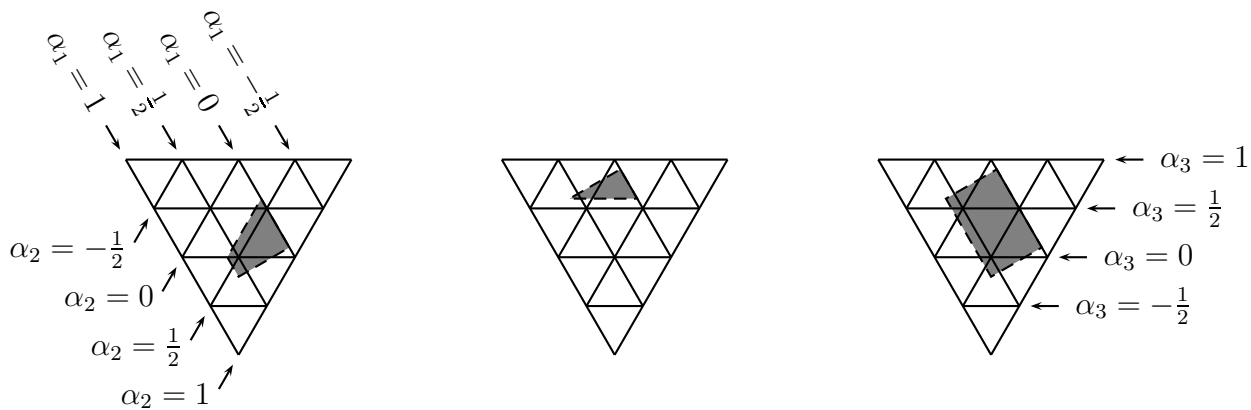
\begin{figure}[htbp]
\psset{unit=.75cm}

\begin{pspicture}(-2,2.7321)(4,-5.1961)

\pspolygon[linestyle=dashed,fillcolor=gray,fillstyle=solid](1.8,-1.7321)(2, -2.0785)(2.9, -1.5588)(2.4,-0.6928)


\multido{\n=0+.5,\nn=4+-.5,\nnn=0+-.8660}{4}{\psline(\n,\nnn)(\nn,\nnn)}
\multido{\n=0+1,\nn=2+.5,\nnn=-3.4641+.8660}{4}{\psline(\n,0)(\nn,\nnn)}
\multido{\n=4+-1,\nn=2+-.5,\nnn=-3.4641+.8660}{4}{\psline(\n,0)(\nn,\nnn)}

\multido{\n=-.35+1.00,\nn=-.1+1.0}{4}{\psline{->}(\n,.6062)(\nn,.1732)}
\multido{\n=.15+.50,\nn=.4+.5,\nnn=-1.4722+-.8660,\nnnn=-1.0392+-.8660}{4}{\psline{->}(\n,\nnn)(\nn,\nnnn)}
\rput[Br](0,-1.732){$\alpha_2=-\frac{1}{2}$}
\rput[Br](.5,-2.5981){$\alpha_2=0$}
\rput[Br](1,-3.4641){$\alpha_2=\frac{1}{2}$}
\rput[Br](1.5,-4.3301){$\alpha_2=1$}
\rput[r]{-60}(-.5,.8660){$\alpha_1=1$}
\rput[r]{-60}(.5,.8660){$\alpha_1=\frac{1}{2}$}
\rput[r]{-60}(1.5,.8660){$\alpha_1=0$}
\rput[r]{-60}(2.5,.8660){$\alpha_1=-\frac{1}{2}$}
\end{pspicture}
\hspace{\stretch{.5}}
\begin{pspicture}(0,2.7321)(4,-5.1961)

\pspolygon[linestyle=dashed,fillcolor=gray,fillstyle=solid](1.2,-0.6928)(2.4,-0.6928)(2.1,-0.1732)


\multido{\n=0+.5,\nn=4+-.5,\nnn=0+-.8660}{4}{\psline(\n,\nnn)(\nn,\nnn)}
\multido{\n=0+1,\nn=2+.5,\nnn=-3.4641+.8660}{4}{\psline(\n,0)(\nn,\nnn)}
\multido{\n=4+-1,\nn=2+-.5,\nnn=-3.4641+.8660}{4}{\psline(\n,0)(\nn,\nnn)}

\end{pspicture}
\hspace{\stretch{.5}}
\begin{pspicture}(0,2.7321)(6,-5.1961)

\pspolygon[linestyle=dashed,fillcolor=gray,fillstyle=solid](2,-2.0785)(2.9,-1.5588)(2.1,-0.1732)(1.2,-0.6928)


\multido{\n=0+.5,\nn=4+-.5,\nnn=0+-.8660}{4}{\psline(\n,\nnn)(\nn,\nnn)}
\multido{\n=0+1,\nn=2+.5,\nnn=-3.4641+.8660}{4}{\psline(\n,0)(\nn,\nnn)}
\multido{\n=4+-1,\nn=2+-.5,\nnn=-3.4641+.8660}{4}{\psline(\n,0)(\nn,\nnn)}
\multido{\n=4.2+-.5,\nn=4.7+-.5,\nnn=0+-.8660}{4}{\psline{<-}(\n,\nnn)(\nn,\nnn)}


\rput[l](5,0){$\alpha_3=1$}
\rput[l](4.5,-.8660){$\alpha_3=\frac{1}{2}$}
\rput[l](4,-1.732){$\alpha_3=0$}
\rput[l](3.5,-2.5981){$\alpha_3=-\frac{1}{2}$}

\end{pspicture}

\caption{Range of exponents  $A^{R}_1$, $A^{R}_2$, $A^{R}$ for $\frac{1}{R}=.4 <\frac{1}{2}$.}\label{ARR>2}
\end{figure}

We first prove the Theorem \ref{VVBHTPmess}. We use results proved in Section \ref{grtLPsec}, with  $\l^P =l^P$. Given sequences of functions $f =\{ f_j\}, g= \{ g_j \}$ define $F(j,x) = f_j(x), G(j,x) =g_j(x)$. Also define the operator $T(F,G)(j,x) = T_j(f_j, g_j)(x)$. Now we need to prove the following bound for $p, q, r, R$ as in the Theorem \ref{VVBHTPmess},
$$ \|T(F, G)\|_{L^r(l^R)} \lesssim \| F \|_{L^p(l^{\infty})} \|G\|_{L^q(l^R)}.$$
Consider the trilinear form associated to $T$,
$$\Lambda(F, G, H) = \int \sum_j T(F,G)(j, x) H(j, x) dx$$
First for $\frac{4}{3} < R < 4$ define the following sets  $A^R_1, A^R_2 \subset \pi_0$,
\begin{equation}\label{AR1}
A^R_1= \left\{ \alpha\in \pi_0: \alpha_1>0, \alpha_2 > \frac{1}{R}, \frac{1}{2} -\frac{1}{2}\alpha_1 -\left|\alpha_2-\frac{1}{2}\right| - \left| \frac{1}{R}-\frac{1}{2}\right| >0, \left| \frac{1}{R} -\frac{1}{2} \right| + \frac{1}{2} \alpha_1 < \frac{1}{R}  \right\}.
\end{equation}
\begin{equation}\label{AR2}
A^R_2= \{\alpha =(\alpha_1, \alpha_2, \alpha_3) \in \pi_0 : (\alpha_1, \alpha_3, \alpha_2) \in A^{R'}_1 \}.
\end{equation}
These  sets are shown in the Figures \ref{ARR=2},  \ref{ARR<2}, and \ref{ARR>2}. Note that  the convex hull of $A^R_1$ and $A^R_2$ is $A^R$.
\begin{equation}\label{AR}
A^R= \text{Convex} (A^R_1 \cup A^R_2).
\end{equation}

We denote $\vec{\alpha}_{R}=(0, \frac{1}{R}, \frac{1}{R'})$. We will show that  $\Lambda $ is of generalized restricted type $(\vec{\alpha}_R, \alpha)$ for $\alpha \in A^R_1$. By doing the same argument with interchanging  second and third coordinates we get $\Lambda$ is of generalized restricted type $(\vec{\alpha}_{R'}, \alpha)$ for $\alpha \in A^{R}_2$.  So to get the Theorem \ref{VVBHTPmess} it is enough to show that $\Lambda $ is generalized restricted type $(\vec{\alpha}_R, \alpha)$ for  $\alpha \in A^R_1$.

Let  $\alpha \in A^R_1$. Note that for $\alpha =(\alpha_1, \alpha_2, \alpha_3) \in A_1^R$ we have $ \alpha_2 - \frac{1}{R}>0$. Given $E=(E_1, E_2, E_3)$  we remove exceptional set from $E_3$,
 $$E_3'=E_3 \setminus \left\{ x: M_{(1+ \epsilon)}1_{E_2} \gtrsim  \left(\frac{|E_2|}{ |E_3|} \right)^{\frac{1}{1+\epsilon}}  \text{or } M_{(1+ \epsilon)}1_{E_1} \gtrsim  \left(\frac{|E_1|}{ |E_3|} \right)^{\frac{1}{1+\epsilon}}  \right\}.
$$
Note that $E'=(E_1, E_2, E_3)$ is a major triple of $(E, \alpha)$. Let $F \in X_{l^P}(E'_1)$ ,$G \in X_{l^Q}(E'_2)$, and $H \in X_{l^R}(E'_3)$ where $X_{l^P}$ as define in \eqref{XLP}. We need to show that,
$$|\Lambda(F, G, H)|= |E|^{\alpha}.$$
Note that we have,
$$\sum_j \|g_j\|_{L^R}\|h_j\|_{L^{R'}} \le \|G\|_{l^R(L^R)} \|H\|_{l^{R'}(L^{R'})} = \|G\|_{L^R(l^R)} \|H\|_{L^{R'}(l^{R'})} \le |E_2|^{\frac{1}{R}} |E_3|^{\frac{1}{R'}}$$
So it is enough to show that,
$$ \int \sum_j |T(F,G)(j, x) H(j, x)| dx \lesssim |E|^{\alpha} |E_2|^{-\frac{1}{R}} |E_3|^{-\frac{1}{R'}} \sum_j \|g_j\|_{L^R}\|h_j\|_{L^{R'}}$$
Now it is enough to prove estimates for single operators,
$$ \int  |T_j(f 1_{E_1},g1_{E_2})( x) h(x) 1_{E_3'}| dx \lesssim |E|^{\alpha} |E_2|^{-\frac{1}{R}} |E_3|^{-\frac{1}{R'}}\|g\|_{L^R}\|h\|_{L^{R'}}.$$
We define the bilinear operators $S_j$ by,
$$S_j(f,g)(x) = T_j(f 1_{E_1},g1_{E_2})( x)1_{E_3'}.$$
Note that it is enough to prove the estimates,
$$\|S_j(f,g)\|_{L^{R}} \lesssim  |E|^{\alpha} |E_2|^{-\frac{1}{R}} |E_3|^{-\frac{1}{R'}} \| f\|_{L^{\infty}} \|g\|_{L^R}$$

We prove this by proving associated trilinear form for $S_j$ is of generalized restricted type $\beta \in \pi_0$ with constant $ |E|^{\alpha} |E_2|^{-\frac{1}{R}} |E_3|^{-\frac{1}{R'}}$ for $\beta $ in a neighborhood of $(0, \frac{1}{R}, \frac{1}{R'})$.

Given a triple $F=(F_1, F_2, F_3)$ we remove exceptional set from $F_1$ and define
 $$F_1'=F_1 \setminus \left\{ x: M_{(1+ \epsilon)}1_{F_2} \gtrsim  \left(\frac{|F_2|}{ |F_1|} \right)^{\frac{1}{1+\epsilon}}  \text{or } M_{(1+ \epsilon)}1_{F_3} \gtrsim  \left(\frac{|F_3|}{ |F_1|} \right)^{\frac{1}{1+\epsilon}}  \right\}
$$
Now as before by using time frequency estimates we get, for $|f| \le 1_{F'_1},$ $ |g| \le 1_{F_2},$ $ |h| \le 1_{F_3}$, and $a=(a_1, a_2, a_3) \in \pi_2$
$$| \lb S_j(f,g), h \rb | \lesssim |F_1|^{a_1} |F_2|^{a_2} |F_3|^{a_3} \sigma_1^{1-2a_1} \sigma_2^{1-2a_2} \sigma_3^{1-2a_3} $$

with following bounds for $\sigma_i$'s,
$$\sigma_1 \lesssim  \left( \frac{|E_1|}{|E_3|} \right)^{\frac{1}{1+\epsilon}} ,$$
$$ \sigma_2 \lesssim \text{min} \left\{ \left( \frac{|E_2|}{|E_3|} \right)^{\frac{1}{1+\epsilon}}, \left( \frac{|F_2|}{|F_1|} \right)^{\frac{1}{1+\epsilon}} \right\},$$
and $$\sigma_3 \lesssim \left( \frac{|F_3|}{|F_1|} \right)^{\frac{1}{1+\epsilon}}.$$
So we have for $\mu, \nu \in [0, \frac{1}{1+\epsilon}]$ and $ \delta_1, \delta_2 >0$ with $ \delta_1 + \delta_2 \in [0,  \frac{1}{1+\epsilon}]$,
$$| \lb S_j(f,g), h \rb | \lesssim |F_1|^{a_1} |F_2|^{a_2} |F_3|^{a_3} \left( \frac{|E_1|}{|E_3|} \right)^{\mu(1-2a_1)}  \left( \frac{|F_2|}{|F_1|} \right)^{\delta_1(1-2a_2)}$$
$$  \left( \frac{|E_2|}{|E_3|} \right)^{\delta_2(1-2a_2)} \left( \frac{|F_3|}{|F_1|} \right)^{\nu(1-2a_3)}.$$
We are need to get estimates of the form,
$$| \lb S_j(f,g), h \rb | \lesssim |E_1|^{\alpha_1} |E_2|^{\alpha_2 - \frac{1}{R}} |E_3|^{\alpha_3 - \frac{1}{R'}} |F_1|^{\beta_1} |F_2|^{\beta_2} |F_3|^{\beta_3},$$
where $\beta =(\beta_1, \beta_2, \beta_3) \in \pi_0$ in a neighborhood of the point $(0, \frac{1}{R}, \frac{1}{R'})$.

First note that the choice of $\mu, \nu, \delta_1 + \delta_2 \in (0, \frac{1}{1+\epsilon})$ allow us reduce problem to choice of $a=(a_1, a_2, a_3) \in \pi_2$ under  the following four conditions. First, third and fourth conditions are obtained from exponents in $|E_1|, |F_3|$ and $|E_2|$. Second condition is obtained by the sum of exponents of $|E_2|$ and $|F_2|$.
\bi
\item $a_1 \in (0, a_1^*) := (0, \frac{1}{2} - \frac{1}{2}\alpha_1 )$
\item $ a_2 \in (0, a_2^*) := (0, \frac{1}{2}- |\alpha_2 -\frac{1}{2} | )$
\item $a_3 \in (0, a_3^*) := (0, \frac{1}{2} - |\frac{1}{R'}-\frac{1}{2}|)$
\item $\alpha_2 - \frac{1}{R} = \delta_2(1-2a_2).$
\ei
First note that we can get $\alpha \in \pi_2 $ satisfying first three conditions if we have,   $a^*_1+ a^*_2 + a^*_3 -1>0.$ This gives us the following condition on $\alpha$ and $ R,$
\begin{equation}\label{maincondition}
\Delta =\Delta(\alpha, R) =  \frac{1}{2} -\frac{1}{2}\alpha_1 -\left|\alpha_2-\frac{1}{2}\right| - \left| \frac{1}{R}-\frac{1}{2}\right| >0.
\end{equation}
Note that this is the third condition in the definition \eqref{AR1} of $A^R_1$.
Now we find additional condition needed to finding $a \in \pi_2$ which also satisfy the last condition. Now note that we have to choose $a_2$ in the range $a_2 \in (a^*_2 -\Delta, a_2^*)$.
Using this with the condition from the sum of exponents in $|F_2|$ and $|E_2|$ we get that to find $a_2, \delta_2$ satisfying last condition we need to have the additional condition,
$$ \alpha_2 -\frac{1}{R} =\delta_2(1-2a_2) < \alpha_2 -(a^*_2 -\Delta)$$
that is,
\begin{equation}\label{maincondition2}
\left| \frac{1}{R} -\frac{1}{2} \right| + \frac{1}{2} \alpha_1 < \frac{1}{R}.
\end{equation}
Note this is the fourth condition on \eqref{AR1}. This completes the proof of the Theorem \ref{VVBHTPmess}.

Next we prove Theorem \ref{VVBHTPmess1}. By using interpolation results proved in Section \ref{grtLPsec} it is enough to show the generalized restricted type condition for associated trilinear form in the corresponding range obtained by interchanging first and second coordinates from the  Theorem \ref{VVBHTPmess}, see \eqref{range}. Note that $BP_j$ is not symmetric in first and second coordinates but we have the same estimates for bilinear operators $S_j$, thus the proof follows the same lines. This completes the proof of the Theorem \ref{VVBHTPmess1}.
\begin{figure}[htbp]
\psset{unit=.75cm}
\begin{pspicture}(-2,2.7321)(4,-5.1961)

\pspolygon[linestyle=dashed,fillcolor=gray,fillstyle=solid](1.5,-2.5980)(2.5,-0.8660)(3,-1.7321)

\multido{\n=0+.5,\nn=4+-.5,\nnn=0+-.8660}{4}{\psline(\n,\nnn)(\nn,\nnn)}
\multido{\n=0+1,\nn=2+.5,\nnn=-3.4641+.8660}{4}{\psline(\n,0)(\nn,\nnn)}
\multido{\n=4+-1,\nn=2+-.5,\nnn=-3.4641+.8660}{4}{\psline(\n,0)(\nn,\nnn)}

\multido{\n=-.35+1.00,\nn=-.1+1.0}{4}{\psline{->}(\n,.6062)(\nn,.1732)}
\multido{\n=.15+.50,\nn=.4+.5,\nnn=-1.4722+-.8660,\nnnn=-1.0392+-.8660}{4}{\psline{->}(\n,\nnn)(\nn,\nnnn)}
\rput[Br](0,-1.732){$\alpha_2=-\frac{1}{2}$}
\rput[Br](.5,-2.5981){$\alpha_2=0$}
\rput[Br](1,-3.4641){$\alpha_2=\frac{1}{2}$}
\rput[Br](1.5,-4.3301){$\alpha_2=1$}
\rput[r]{-60}(-.5,.8660){$\alpha_1=1$}
\rput[r]{-60}(.5,.8660){$\alpha_1=\frac{1}{2}$}
\rput[r]{-60}(1.5,.8660){$\alpha_1=0$}
\rput[r]{-60}(2.5,.8660){$\alpha_1=-\frac{1}{2}$}
\end{pspicture}
\hspace{\stretch{.5}}
\begin{pspicture}(0,2.7321)(4,-5.1961)

\pspolygon[linestyle=dashed,fillcolor=gray,fillstyle=solid](.5,-0.8660)(2,0)(2.5,-0.8660)

\multido{\n=0+.5,\nn=4+-.5,\nnn=0+-.8660}{4}{\psline(\n,\nnn)(\nn,\nnn)}
\multido{\n=0+1,\nn=2+.5,\nnn=-3.4641+.8660}{4}{\psline(\n,0)(\nn,\nnn)}
\multido{\n=4+-1,\nn=2+-.5,\nnn=-3.4641+.8660}{4}{\psline(\n,0)(\nn,\nnn)}

\end{pspicture}
\hspace{\stretch{.5}}
\begin{pspicture}(0,2.7321)(6,-5.1961)

\pspolygon[linestyle=dashed,fillcolor=gray,fillstyle=solid](.5,-0.8660)(2,0)(3,-1.7321)(1.5,-2.5980)


\multido{\n=0+.5,\nn=4+-.5,\nnn=0+-.8660}{4}{\psline(\n,\nnn)(\nn,\nnn)}
\multido{\n=0+1,\nn=2+.5,\nnn=-3.4641+.8660}{4}{\psline(\n,0)(\nn,\nnn)}
\multido{\n=4+-1,\nn=2+-.5,\nnn=-3.4641+.8660}{4}{\psline(\n,0)(\nn,\nnn)}
\multido{\n=4.2+-.5,\nn=4.7+-.5,\nnn=0+-.8660}{4}{\psline{<-}(\n,\nnn)(\nn,\nnn)}


\rput[l](5,0){$\alpha_3=1$}
\rput[l](4.5,-.8660){$\alpha_3=\frac{1}{2}$}
\rput[l](4,-1.732){$\alpha_3=0$}
\rput[l](3.5,-2.5981){$\alpha_3=-\frac{1}{2}$}

\end{pspicture}

\caption{Range of exponents  $A^{R}_1$, $A^{R}_2$, $A^{R}$ for $R=2$.}\label{ARR=2}
\end{figure}
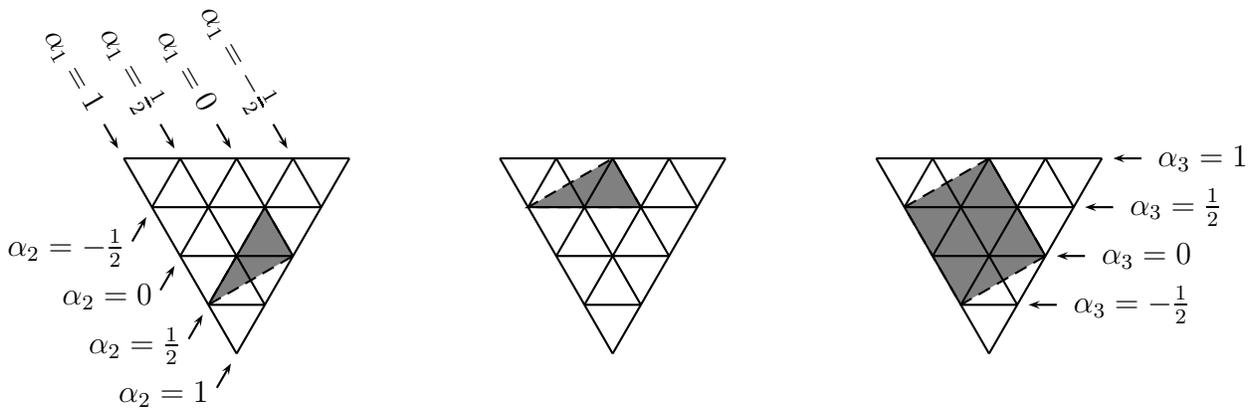



\newpage
\begin{thebibliography}{03}

\bibitem{BT1} Bateman M., Thiele C.,
\emph{ $L^p$ estimates for the Hilbert transform along a one-variable vector field}
preprint vailable on arxiv

\bibitem{BP1} A. Benedek and R. Panzone,
\emph{ The space $L^P$, with mixed norm}
Duke Math. J. Volume 28, Number 3 (1961), 301-324.




\bibitem{GL} Grafakos L. and Li X.,
\emph{ The Disc as a Bilinear Multiplier}
American Journal of Mathematics
Vol. 128, No. 1 (Feb., 2006), pp. 91-119

\bibitem{LL1} Lacey M. and  Li X.,
\emph{  Maximal Theorems for The Directional Hilbert transform}
Transactions of the American Mathematical Society
Volume 358, Number 9, September 2006, Pages 4099-4117


\bibitem{LT1} Lacey M. and Thiele C.,
\emph{  On Calderon's Conjecture}
Annals of Mathematics, 149(1999), 475-496

\bibitem{MPTT} Muscalu C., Pipher J.,  Tao T. and Thiele C.,
\emph{  Bi-parameter paraproducts}
Acta Mathematica Volume 193, Number 2, 269-296

\bibitem{MTT1} Muscalu C., Tao T. and Thiele C.,
\emph{  Multi-linear operators given by singular multipliers}
J. Amer. Math. Soc. 15 (2002), 469-496.

\bibitem{MTT2} Muscalu C., Tao T. and Thiele C.,
\emph{  $L^p$ estimates for the biest II. The Fourier case}
Math. Ann. 329 (2004), no. 3, 427-461

\bibitem{Th1} Thiele C.,
\emph{ Wave Packet Analysis (CBMS Regional Conference Series in Mathematics)}



\end{thebibliography}
\end{document}